\documentclass[12pt]{article}%
\pdfoutput = 1
\usepackage[left=2.2cm, right=2.2cm, top=2.2cm, bottom=2.2cm]{geometry}
\usepackage{amsfonts}
\usepackage{amssymb}
\usepackage{amsmath}
\usepackage{lscape}%
\usepackage{graphicx}
\providecommand{\U}[1]{\protect\rule{.1in}{.1in}}

\begin{document}

\title{Approximate Bayesian Computation in State Space Models\thanks{This research
has been supported by Australian Research Council (ARC) Future Fellowship
FT0991045. We appreciate the input of various participants at \textit{ABC in
Rome}, June 2013, and \textit{ABC in Sydney}, July 2014, on earlier drafts of
the paper.}}
\author{Gael M. Martin\thanks{Department of Econometrics and Business Statistics,
Monash University, Melbourne, Australia. Corresponding author; email:
gael.martin@monash.edu.}, Brendan P.M. McCabe\thanks{Management School,
University of Liverpool, U.K.}, Worapree Maneesoonthorn\thanks{Melbourne
Business School, University of Melbourne, Australia.}
\and and Christian P. Robert\thanks{University of Paris Dauphine, Centre de
Recherche en \'{E}conomie et Statistique, and University of Warwick.}}
\maketitle

\begin{abstract}
A new approach to inference in state space models is proposed, based on
approximate Bayesian computation (ABC). ABC avoids evaluation of the
likelihood function by matching observed summary statistics with statistics
computed from data simulated from the true process; exact inference being
feasible only if the statistics are sufficient. With finite sample sufficiency
unattainable in the state space setting, we seek asymptotic sufficiency via
the maximum likelihood estimator (MLE) of the parameters of an auxiliary
model. We prove that this auxiliary model-based approach achieves Bayesian
consistency, and that - in a precise limiting sense - the proximity to
(asymptotic) sufficiency yielded by the MLE is replicated by the score. In
multiple parameter settings a separate treatment of scalar parameters, based
on integrated likelihood techniques, is advocated as a way of avoiding the
curse of dimensionality. Some attention is given to a structure in which the
state variable is driven by a continuous time process, with exact inference
typically infeasible in this case as a result of intractable transitions. The
ABC method is demonstrated using the unscented Kalman filter as a fast and
simple way of producing an approximation in this setting, with a stochastic
volatility model for financial returns used for illustration.\medskip

\noindent\emph{Keywords:} Likelihood-free methods, latent diffusion models,
linear Gaussian state space models, asymptotic sufficiency, unscented Kalman
filter, stochastic volatility.

\medskip

\noindent\emph{JEL Classification:} C11, C22, C58

\end{abstract}

\newpage

\section{Introduction}

\baselineskip18pt

Approximate Bayesian computation (ABC) (or likelihood-free inference) has
become increasingly prevalent in areas of the natural sciences in which
likelihood functions requiring integration over a large number of complex
latent states are essentially intractable. (See Cornuet \textit{et al., }2008,
Beaumont, 2010, Marin \emph{et al.}, 2011 and Sisson and Fan, 2011 for recent
reviews.) The technique circumvents direct evaluation of the likelihood
function by matching summary statistics calculated from the observed data with
corresponding statistics computed from data simulated from the assumed data
generating process. If such statistics are sufficient, the method yields an
approximation to the exact posterior distribution of interest that is
accurate, given an adequate number of simulations; otherwise, \textit{partial
}posterior inference, reflecting the information content of the set of summary
statistics only, is the outcome.

The choice of statistics for use within the ABC method, in addition to
techniques for determining the matching criterion, are clearly of paramount
importance, with much recent research having been devoted to devising ways of
ensuring that the information content of the chosen set is maximized, in some
sense; e.g.\textbf{\ }Joyce and Marjoram (2008), Wegmann \textit{et al.}
(2009), Blum (2010a) and Fearnhead and Prangle (2012). Recent contributions
here include those of Drovandi \textit{et al.} (2011), Drovandi and Pettitt
(2013), Gleim and Pigorsch (2013) and Creel and Kristensen (2014), in which
the statistics are produced by estimating an approximating \textit{auxiliary
}model using both simulated and observed data. This approach mimics, in a
Bayesian framework, the principle underlying the frequentist methods of
indirect inference (II) (Gouri\'{e}roux \textit{et al.} 1993, Smith, 1993,
Heggland and Frigessi, 2004) and efficient method of moments (EMM) (Gallant
and Tauchen, 1996), using, as it does, the approximating model to produce
feasible, but sub-optimal, inference about an intractable true model. Whilst
the price paid for the approximation in the frequentist setting is a possible
reduction in efficiency, the price paid in the Bayesian case is posterior
inference that is conditioned on statistics that are not sufficient for the
parameters of the true model, and which amounts to only partial inference as a consequence.

Our paper continues in this spirit, but with particular focus given to the
application of auxiliary model-based ABC methods in the state space model
(SSM) framework. We begin by demonstrating that reduction to a set of
sufficient statistics of fixed dimension relative to the sample size is
\textit{infeasible} in finite samples in SSMs. This key observation then
motivates our decision to seek asymptotic sufficiency in the state space
setting by using the MLE of the parameters of the auxiliary model as the
(vector) summary statistic in the ABC matching criterion. We focus on two
qualitatively different cases: 1) one in which the auxiliary model
\textit{coincides} with the true model, in which case asymptotic sufficiency
for the true parameters is achievable via the proposed ABC technique; and 2)
the more typical case in which the exact likelihood function is inaccessible,
and the auxiliary model represents an approximation only. The first case
mimics that sometimes referenced in the II (or EMM) literature, in which the
auxiliary model `nests', or is equivalent to in some well-defined sense, the
true model, and full asymptotic efficiency is achieved by the frequentist
methods as a consequence. Investigation of this case allows us to document the
maximum accuracy gains that are possible via the auxiliary model route,
compared with ABC techniques based on alternative summaries, without the
confounding effect of the error in the approximating model. The second case
gives some insight into what can be achieved in a general non-linear state
space setting when the investigator is forced to adopt an inexact
approximating model in the implementation of auxiliary model-based ABC. We
give emphasis here to non-linear models in which the state (and possibly the
observed) is driven by a continuous time model, as this is the canonical
example in which simulation from the true model is feasible (at least via an
arbitrarily fine discretization), whilst the likelihood function is
(typically) unavailable and exact posterior analysis thus not achievable.

We begin by considering the very concept of finite sample sufficiency in the
state space context, and the usefulness of applying a typical ABC approach -
based on \textit{ad hoc} summary statistics - in this setting. Using the
linear Gaussian model for illustration, we demonstrate the lack of reduction
to a set of sufficient statistics of fixed dimension, this result providing
motivation, as noted above, for the pursuit of asymptotic sufficiency via the
auxiliary model method. We then proceed to demonstrate the Bayesian
consistency of the auxiliary model approach, subject to the typical quasi-MLE
form of conditions being satisfied. We also illustrate that to the order of
accuracy that is relevant in establishing the theoretical properties of an ABC
technique (i.e. allowing the tolerance used in the matching of the statistics
to approach zero), a selection criterion based on the score of the auxiliary
model - evaluated at the MLE computed from the observed data - yields
equivalent results to a criterion based directly on the MLE itself. This
equivalence is shown to hold in both the exactly and over-identified cases,
and independently of any (positive definite) weighting matrix used to define
the two alternative distance measures, and implies that the proximity to
asymptotic sufficiency yielded by the use of the MLE in an ABC algorithm will
be replicated by the use of the score. Given the enormous gain in speed
achieved by avoiding optimization of the approximate likelihood at each
replication of ABC, this is an critical result from a computational
perspective. The application of the proposed method in multiple parameter
settings is addressed, with separate treatment of scalar (or lower-dimensional
blocks of) parameters, via marginal, or integrated likelihood principles
advocated, as a possible way of avoiding the inaccuracy that plagues ABC
techniques in high dimensions. (See Blum, 2010b and Nott \emph{et al.}, 2014).

The results outlined in the previous paragraph are applicable to an auxiliary
model-based ABC method applied in any context (subject to regularity) and,
hence, are of interest in their own right. However, our particular interest,
as already noted, is in applying the auxiliary model method - and thereby
exploiting these properties \ - in the state space setting. For Case 1) we
choose to illustrate the approach using the linear Gaussian model, whereby the
exact likelihood (and hence score) is accessible via the Kalman filter (KF),
and asymptotic sufficiency thus achievable. For Case 2) we illustrate the
approach via a particular choice of (approximating) auxiliary model for the
continuous time SSM. Specifically, the approximating model is formed as a
discretization of the true continuous time model, with the augmented unscented
Kalman filter (AUKF) (Julier \textit{et al.}, 1995, Julier and Uhlmann,
2004)\ used to evaluate the likelihood of that model. The general
applicability, speed and simplicity of the AUKF calculations render the ABC
scheme computationally feasible and relatively simple to implement. This
particular approach to the definition of an auxiliary model also leads to a
set of summary statistics of relatively small dimension. This is in contrast,
for example, with an approach based on a highly parameterized (`nesting')
approximating model (see, for example, Gleim and Pigorsch, 2013), in which the
large number of auxiliary parameters - in principle sufficient for the
parameters of the true latent diffusion model - is likely to yield a very
inaccurate (non-parametric) estimate of the true posterior, due to the large
dimension of the conditioning statistics. The equality between the number of
parameters in our exact and approximating models also means that
marginalization of the approximating model to produce a scalar matching
criterion for each parameter of the true model is meaningful.

The paper proceeds as follows. In Section \ref{abc} we briefly summarize the
basic principles of ABC as they would apply in a state space setting,
including the role played by summary statistics and sufficiency. We
demonstrate the lack of finite sample sufficiency reduction in an SSM, using
the linear Gaussian model for illustration.\textbf{\ }In Section \ref{Aux}, we
then proceed to demonstrate the properties of ABC based on the MLE, score and
marginal score, respectively, of a generic approximating model followed, in
Section \ref{model}, by an outline of a computationally feasible approximation
- based on the AUKF - for use in the non-linear state space setting. Using
(repeated samples of) artificially generated data, the accuracy with which the
proposed technique reproduces the exact posterior distribution is assessed in
Section \ref{assess}. The ABC methods are based respectively on: i) the joint
score; ii) the marginal score; iii) a (weighted) Euclidean metric based on
statistics that are sufficient for an observed autoregressive model of order
one; and iv) the dimension-reduction technique of Fearnhead and Prangle
(2012), applied to this latter set of summary statistics. We conduct the
assessment firstly within the context of the linear Gaussian model, with the
issues of sufficiency and matching that are key to accurately reproducing the
true posterior distribution able to be illustrated precisely in this setting.
The overall superiority of both the joint and marginal score techniques over
the ABC methods based on summary statistics is demonstrated numerically, as is
the remarkable accuracy yielded by the marginal score technique in particular.
This exercise thus forms a resounding proof-of-concept for the score-based ABC
method, albeit in a case where the exact score is available.

We then proceed to assess performance in a particular non-linear latent
diffusion model, in which the degree of accuracy of the AUKF-based
approximating model plays a role. A stochastic volatility for financial
returns, in which the latent volatility is driven by a square root diffusion
model, is adopted as the non-linear example, as the existence of known
(non-central chi-squared) transition densities means that the exact likelihood
function/posterior distribution is available, for the purpose of comparison.
We apply the deterministic grid-based filtering method of Ng \textit{et al.}
(2013) - suitable for this particular setting - to produce the exact
comparators for our ABC-based estimates of the relevant marginal posteriors,
as well as the marginal posteriors associated with an Euler approximation to
the true model. The score methods out-perform the summary statistic methods in
the great majority of cases documented. Some gain in accuracy is still
produced via the marginalization technique, although that gain is certainly
less marked than in the linear Gaussian case, in which the exact score is
accessible. Notably, all ABC-based approximations, which exploit simulation
from the exact latent diffusion model, serve as more accurate estimates
(overall) of the exact posteriors than do the AUKF and Euler approximations
themselves. Section \ref{end} concludes.

\section{ABC in State Space Models\label{abc}}

\subsection{Outline of the basic approach}

The aim of ABC is to produce draws from an approximation to the posterior
distribution of a vector of unknowns, $\mathbf{\theta}$, given the
$T$-dimensional vector of observed data $\mathbf{y}=(y_{1},y_{2}%
,...,y_{T})^{\prime}$,
\[
p(\mathbf{\theta|y})\propto p(\mathbf{y|\theta})p(\mathbf{\theta}),
\]
in the case where both the prior, $p(\mathbf{\theta})$, and the likelihood,
$p(\mathbf{y|\theta})$, can be simulated. These draws are used, in turn, to
approximate posterior quantities of interest, including marginal posterior
moments, marginal posterior distributions and predictive distributions. The
simplest (accept/reject) form of the algorithm (Tavar\'{e}\textit{\ et al.
}1997, Pritchard, 1999) proceeds as follows:

\begin{enumerate}
\item[Step 1] Simulate $\mathbf{\theta}^{i}$, $i=1,2,...,R$, from
$p(\mathbf{\theta})$

\item[Step 2] Simulate $\mathbf{z}^{i}=(z_{1}^{i},z_{2}^{i},...,z_{T}%
^{i})^{\prime}$, $i=1,2,...,R$, from the likelihood, $p(\mathbf{.|\theta}%
^{i})$

\item[Step 3] Select $\mathbf{\theta}^{i}$ such that:%
\begin{equation}
d\{\mathbf{\eta}(\mathbf{y}),\mathbf{\eta}(\mathbf{z}^{i})\}\leq
\varepsilon,\label{distance}%
\end{equation}
where $\eta(\mathbf{.})$ is a (vector) statistic, $d\{.\}$ is a distance
criterion of some sort, and the tolerance level $\varepsilon$ is arbitrarily
small. In practice $\varepsilon$ may be chosen such that, for a given value of
$R$, a certain (small) proportion of draws of $\mathbf{\theta}^{i}$ are selected.
\end{enumerate}

The algorithm thus samples $\mathbf{\theta}$ and $\mathbf{z}$ from the
joint\textit{\ }posterior:%
\[
p_{\varepsilon}(\mathbf{\theta},\mathbf{z|y})=\frac{p(\mathbf{\theta
})p(\mathbf{z|\theta})\mathbb{I}_{A_{\varepsilon,\mathbf{y}}}(\mathbf{z})}{%
%TCIMACRO{\tint }%
%BeginExpansion
{\textstyle\int}
%EndExpansion
p(\mathbf{\theta})p(\mathbf{z|\theta})\mathbb{I}_{A_{\varepsilon,\mathbf{y}}%
}(\mathbf{z})},
\]
where $\mathbb{I}_{B}$ is the indicator function defined on the set $B$ and
$A_{\varepsilon,\mathbf{y}}=\{\mathbf{z}\in\mathbf{S};d\{\mathbf{\eta
}(\mathbf{y}),\eta(\mathbf{z})\}\leq\varepsilon\}.$ Clearly, when
$\eta(\mathbf{.})$ is sufficient and $\varepsilon$ arbitrarily small,%
\[
p_{\varepsilon}(\mathbf{\theta|y})=%
%TCIMACRO{\tint }%
%BeginExpansion
{\textstyle\int}
%EndExpansion
p_{\varepsilon}(\mathbf{\theta},\mathbf{z|y})d\mathbf{z}%
\]
approximates the true posterior, $p(\mathbf{\theta|y})$, and draws from
$p_{\varepsilon}(\mathbf{\theta},\mathbf{z|y})$ can be used to estimate
features of the true posterior. In practice however, the complexity of the
models to which ABC is applied implies, almost by definition, that sufficiency
is unattainable. Hence, in the limit, as $\varepsilon\rightarrow0$, the draws
can be used only to approximate features of $p(\mathbf{\theta|\eta}%
(\mathbf{y})).$

Adaptations of the basic rejection scheme have involved post-sampling
corrections of the draws using kernel methods (Beaumont \textit{et al.}, 2002,
Blum, 2010a, Blum and Fran\c{c}ois, 2010), or the insertion of Markov chain
Monte Carlo (MCMC) and/or sequential Monte Carlo (SMC) steps (Marjoram
\textit{et al.}, 2003, Sisson \textit{et al.,} 2007, Beaumont \textit{\ et
al.}, 2009, Toni \textit{et al.}, 2009 and Wegmann \textit{et al.}, 2009), to
improve the accuracy with which $p(\mathbf{\theta|\eta}(\mathbf{y}))$ is
estimated, for any given number of draws. Focus is also given to choosing
$\eta(.)$ and/or $d\{.\}$ so as to render $p(\mathbf{\theta|\eta}%
(\mathbf{y}))$ a closer match to $p(\mathbf{\theta|y})$, in some sense; see
Joyce and Marjoram (2008), Wegmann \textit{et al.}, Blum (2010b) and Fearnhead
and Prangle (2012). In the latter vein, Drovandi \textit{et al.} (2011) argue,
in the context of a specific biological model, that the use of $\eta(.)$
comprised of the MLEs of the parameters of a well-chosen approximating model,
may yield posterior inference that is conditioned on a large portion of the
information in the data and, hence, be close to exact inference based on
$p(\mathbf{\theta|y})$. (See also Drovandi and Pettitt, 2013, Gleim and
Pigorsch, 2013, and Creel and Kristensen, 2014). It is the spirit of this
approach that informs the current paper, but with our attention given to
rendering the approach feasible in a \textit{general} state space framework
that encompasses a large number of the models that are of interest to
practitioners, including continuous time models.

Our focus then is on the application of ABC in the context of a general SSM
with measurement and transition distributions,%
\begin{align}
& p(y_{t}|x_{t},\mathbf{\phi})\label{meas_gen}\\
& p(x_{t}|x_{t-1},\mathbf{\phi})\label{state_gen}%
\end{align}
respectively, where $\mathbf{\phi}$ is a $p$-dimensional vector of static
parameters, elements of which may characterize either the measurement or state
relation, or both.\emph{\ }For expositional simplicity, and without loss of
generality, we consider the case where both $y_{t}$ and $x_{t}$ are scalars.
In financial applications it is common that both the observed and latent
processes\emph{\ }are driven by continuous time processes, with the transition
distribution in (\ref{state_gen}) being unknown (or, at least, computationally
challenging) as a consequence. Bayesian inference would then typically proceed
by invoking (Euler) discretizations for both the measurement and state
processes and applying MCMC- or SMC-based techniques, with such methods being
tailor-made to suit the features of the particular (discretized) model at
hand; see Giordini \textit{et al.} (2011) for a recent review. In some models
expressed initially in discrete time, it may also be the case that the
conditional distribution in (\ref{meas_gen}) is unavailable in closed form,
such as when empirically relevant distributions for financial returns are
adopted (e.g. Peters \emph{et al.,} 2012). In such cases, SMC- or MCMC-based
inferential methods are typically infeasible. In contrast, in all of these
cases the proposed ABC method \textit{is} feasible, as long as simulation from
the true model (at least via an arbitrarily fine discretization, in the
continuous time case) is possible.

The full set of unknowns thus constitutes the augmented vector $\mathbf{\theta
}=(\mathbf{\phi}^{\prime}\mathbf{,x}_{c}^{\prime})^{\prime}$ where, in the
case where $x_{t}$ evolves in continuous time, $\mathbf{x}_{c}$ represents the
infinite-dimensional vector comprising the continuum of unobserved states over
the sample period. However, to fix ideas, we define $\mathbf{\theta
}=(\mathbf{\phi}^{\prime}\mathbf{,x}^{\prime})^{\prime},$ where $\mathbf{x}%
=(x_{1},x_{2},...,x_{T})^{\prime}$ is the $T$-dimensional vector comprising
the time $t$ states for the $T$ observation periods in the
sample.\footnote{For example, in a continuous time stochastic volatility model
such values may be interpreted as end-of-day volatilities.} Implementation of
the algorithm thus involves simulating from $p(\mathbf{\theta})$ by simulating
$\mathbf{\phi}$ from the prior $p(\mathbf{\phi})\,$, followed by simulation of
$x_{t}$ via the process for the state, conditional on the draw of
$\mathbf{\phi}$, and subsequent simulation of artificial data $z_{t}$
conditional on the draws of $\mathbf{\phi}$ and the state variable.

Crucially, the focus in this paper is on inference about $\mathbf{\phi}%
$\ only; hence, only draws of $\mathbf{\phi}$\ are retained (via the selection
criterion) and those draws used to produce an estimate of the marginal
posterior, $p(\mathbf{\phi|y})$, and with sufficiency to be viewed as relating
to $\mathbf{\phi}$\ only. Hence, from this point onwards, when we reference a
vector of summary statistics, $\mathbf{\eta(y)}$, it is the information
content of that vector with respect to $\mathbf{\phi}$\ that is of importance,
and the proximity of $p(\mathbf{\phi|\eta(y)})$\ to the marginal posterior of
$\mathbf{\phi}$ that is under question. We comment briefly on state inference
in Section \ref{end}.

Before outlining the proposed methodology for the model in (\ref{meas_gen})
and (\ref{state_gen}) in Section \ref{Aux} we highlight the key observation
that motivates our approach, namely that reduction to sufficiency in finite
samples is not possible in state space settings.\textbf{\ }We use a linear
Gaussian state space model to illustrate this result, as closed-form
expressions are available in this case; however, as highlighted at the end of
the section, the result is, in principle, applicable to any SSM.

\subsection{Lack of finite sample sufficiency reduction\label{lack}}

Sufficient statistics are useful for inference about a (vector) parameter
since, being in possession of the sufficient set means that the data itself
may be discarded for inference purposes. When the cardinality of the
sufficient set is small relative to the sample size a significant reduction in
complexity is achieved and in the case of ABC, conditioning on the sufficient
statistics leads to no loss of information, and the method produces a
simulation-based estimate of the true posterior. The difficulty that arises is
that only distributions that are members of the exponential family (EF)
possess sufficient statistics that achieve a reduction to a fixed dimension
relative to the sample size. In the context of the general SSM described by
(\ref{meas_gen}) and (\ref{state_gen}) the effective use of sufficient
statistics is problematic. For any $t$ it is unlikely that the marginal
distribution of $y_{t}$ will be a member of the EF, due to the vast array of
non-linearities that are possible, in either the measurement or state
equations, or both. Moreover, even if $y_{t}$ \textit{were} a member of the EF
for each $t$, to achieve a sufficiency reduction it is required that the
\textit{joint} distribution of $\mathbf{y}=\left\{  y_{t};t=1,2,...,T\right\}
$ also be in the EF. For example, even if $y_{t}$ were Gaussian, it does not
necessarily follow that the joint distribution of $\mathbf{y}$ will achieve a
sufficiency reduction. The most familiar example of this is when $\mathbf{y}$
follows a Gaussian moving average (MA) process and consequently only the whole
sample is sufficient.

Even the simplest SSMs generate MA-like dependence in the data. Consider the
linear Gaussian SSM, expressed in regression form as
\begin{align}
y_{t}  & =x_{t}+e_{t}\label{measlg}\\
x_{t}  & =\delta+\rho x_{t-1}+v_{t},\label{statelg}%
\end{align}
where the disturbances are respectively independent $N\left(  0,\sigma_{e}%
^{2}=1\right)  $ and $N\left(  0,\sigma_{v}^{2}\right)  $ variables. In this
case, the joint distribution of the vector of $y_{t}$'s (which are marginally
normal and members of the EF) is $\mathbf{y\sim N[\mu,}\sigma_{x}^{2}\left(
r\mathbf{I}+\mathbf{V}\right)  ),$ where $r=\sigma_{e}^{2}/\sigma_{x}^{2}$ is
the inverse of the signal-to-noise (SN) ratio and $\mathbf{V}$ is the familiar
Toeplitz matrix associated with an autoregressive (AR) model of order 1. To
construct the sufficient statistics we need to evaluate $\left(
r\mathbf{I}+\mathbf{V}\right)  ^{-1}$, which appears in the quadratic form of
the multivariate normal density, with the structure of $\left(  r\mathbf{I}%
+\mathbf{V}\right)  ^{-1}$ determining the way in which sample information
about the parameters is accumulated and, hence, the sufficiency reduction that
is achievable. (See, for example, Anderson, 1958, Chp 6.) Representing
$\left(  r\mathbf{I}+\mathbf{V}\right)  ^{-1}$ as
\begin{equation}
\left(  r\mathbf{I}+\mathbf{V}\right)  ^{-1}=\mathbf{V}^{-1}-r\mathbf{V}%
^{-2}+r^{2}\mathbf{V}^{-3}-...,\label{expansion}%
\end{equation}
it is straightforward to show that as the order of the approximation is
increased by retaining more terms, the extent of the accumulation across
successive observations is reduced, with the full sample of observations on
$y_{t}$ ultimately being needed to attain sufficiency. Given that the
magnitude of $r$ determines how many terms in (\ref{expansion}) are required
for the approximation to be accurate, we see that the SN ratio determines how
well the set of summary statistics \textit{that would be sufficient} for an
observed AR(1) process (with $r=0$), namely,%

\begin{equation}
s_{1}=\sum_{t=2}^{T-1}y_{t},\text{ }s_{2}=\sum_{t=2}^{T-1}y_{t}^{2},\text{
}s_{3}=\sum_{t=2}^{T}y_{t}y_{t-1},\text{ }s_{4}=y_{1}+\text{ }y_{T},\text{
}s_{5}=y_{1}^{2}+y_{T},\label{AR1_summ_stats}%
\end{equation}
approximates the information content of the true set of sufficient statistics,
that is, the full sample. If the SN ratio is large (i.e. $r$ is small) then
using the set in (\ref{AR1_summ_stats}) as summary statistics may produce a
reasonable approximation to sufficiency. However, as the SN ratio declines
(and higher powers of $r$ cannot be ignored as a consequence) then this set
deviates further and further from sufficiency.

This same qualitative problem would also characterize any SSM nested in
(\ref{meas_gen}) and (\ref{state_gen}), with the only difference being that,
in any particular case there would not necessarily be an analytical link
between the SN ratio and the lack of sufficiency associated with any finite
set of statistics calculated from the observations. The quest for an accurate
ABC technique in a state space setting - based on an arbitrary set of
statistics - is thus not well-founded and this, in turn, motivates the search
for asymptotic sufficiency via the MLE.

\section{Auxiliary model-based ABC \label{Aux}}

\subsection{Theoretical properties\label{Theory}}

The asymptotic Gaussianity of the MLE for the parameters of (\ref{meas_gen})
and (\ref{state_gen}) (under regularity) implies that the MLE\ satisfies the
factorization theorem and is thereby asymptotically sufficient for the
parameters of that model. (See Cox and Hinkley, 1974, Chp. 9 for elucidation
of this matter.) Denoting the log-likelihood function by\textbf{\ }%
$L(\mathbf{y;\phi})$, maximizing\textbf{\ }$L(\mathbf{y;\phi})$ with respect
to $\mathbf{\phi}$ yields $\widehat{\mathbf{\phi}}$, which could, in
principle, be used to define $\mathbf{\eta}(\mathbf{.})$ in an ABC algorithm.
For large enough $T$ the algorithm would produce draws from the exact
posterior. Indeed, in arguments that mirror those adopted by Gallant and
Tauchen (1996) and Gouri\'{e}roux \textit{et al.} (1993) for the EMM and II
estimators respectively, Gleim and Pigorsch (2013) demonstrate that if
$\mathbf{\eta}(\mathbf{.})$ is chosen to be the MLE of an auxiliary model that
`nests' the true model in some well-defined way, asymptotic sufficiency will
still be achieved; see also Gouri\'{e}roux and Monfort (1995) on this point.

Of course, if the SSM in question is such that the exact likelihood is
accessible, the model is likely to be tractable enough to preclude the need
for treatment via ABC. Further, as we allude to in the Introduction, the quest
for asymptotic sufficiency via a (possibly large) nesting auxiliary model
conflicts with the quest for an accurate non-parametric estimate of the
posterior using the ABC draws; a point that, to our knowledge, has not been
noted in the literature. Hence, in practice, the appropriate goal in the ABC
context is to define a \textit{parsimonious}, analytically tractable (and
computationally efficient) model that \textit{approximates} the (generally
intractable) data generating process in (\ref{meas_gen}) and (\ref{state_gen})
as well as possible, and use that model as the basis for constructing a
summary statistic within an ABC algorithm. If the approximating model is
`accurate enough' as a representation of the true model, such an approach will
yield, via the ABC algorithm, an estimate of the posterior distribution that
is conditioned on a statistic that is `close to' being sufficient, at least
for a large enough sample.

Despite the loss of full (asymptotic) sufficiency associated with the use of
an approximating model to generate the matching statistics in an ABC
algorithm, we show here that Bayesian consistency will still be achieved,
subject to certain regularity conditions. Define the choice criterion in
(\ref{distance}) as%

\begin{equation}
d\{\mathbf{\eta}(\mathbf{y}),\mathbf{\eta}(\mathbf{z}^{i})\}=\sqrt{\left[
\widehat{\mathbf{\beta}}(\mathbf{y)-}\widehat{\mathbf{\beta}}(\mathbf{z}%
^{i}\mathbf{)}\right]  ^{\prime}\mathbf{\Omega}\left[  \widehat{\mathbf{\beta
}}(\mathbf{y)-}\widehat{\mathbf{\beta}}(\mathbf{z}^{i}\mathbf{)}\right]  }%
\leq\varepsilon,\label{dist_mle}%
\end{equation}
where $\widehat{\mathbf{\beta}}(\mathbf{.)}$ is the MLE of the parameter
vector $\mathbf{\beta}$ of the auxiliary model with log-likelihood function,
$L_{a}(\mathbf{y;\beta})$, and $\mathbf{\Omega}$ is some positive definite
matrix. The quadratic function under the square root essentially mimics the
criterion used in the II technique, in which case $\mathbf{\Omega}$ would
assume the sandwich form of variance-covariance estimator - appropriate for
when the auxiliary model does not coincide with the true model - and
optimization with respect to the parameters of the latter is the
goal.\footnote{In practice the implementation of II involves the use of a
simulated sample in the computation of $\widehat{\mathbf{\beta}}%
(\mathbf{z}^{i}\mathbf{)}$ that is a multiple of the size of the empirical
sample.} In Bayesian analyses, in which (\ref{dist_mle}) is used to produce
ABC draws, $\mathbf{\Omega}$ may also be defined as the sandwich estimator
(Drovandi and Pettit, 2013, and Gleim and Pigorsch, 2013), or simply as the
inverse of the (estimated) variance-covariance matrix of
$\widehat{\mathbf{\beta}}$, evaluated at $\widehat{\mathbf{\beta}}%
(\mathbf{y)}$ (Drovandi\textit{\ et al.,} 2011). However, in common with the
frequentist proof of consistency of the II estimator, Bayesian consistency -
whereby the posterior for $\mathbf{\phi}$ is degenerate as $T\rightarrow
\infty$ at the true $\mathbf{\phi=\phi}_{0}$ - is invariant to the choice of
the (positive definite) $\mathbf{\Omega.}$ The demonstration follows directly
from the same arguments used to prove consistency of the II estimator (see,
for e.g. Gouri\'{e}roux and Monfort, 1996, Appendix 4A.1), with the following
regularity conditions required:

\begin{description}
\item[(A1)] lim$_{T\rightarrow\infty}T^{-1}L_{a}(\mathbf{z(\phi);\beta
})=L_{\infty}(\mathbf{\phi,\beta})$, uniformly in $\mathbf{\beta}$, where
$L_{\infty}(\mathbf{\phi,\beta})$ is a deterministic limit function.

\item[(A2)] $L_{\infty}(\mathbf{\phi,\beta})$ has a unique maximum with
respect to $\mathbf{\beta:}$ $\mathbf{b}(\mathbf{\phi)=\arg}$
$\underset{\mathbf{\beta}}{\text{max }}L_{\infty}(\mathbf{\phi,\beta}).$

\item[(A3)] The equation $\mathbf{\beta=b(\phi)}$ admits a unique solution in
$\mathbf{\phi}$, for all $\mathbf{\phi}$.
\end{description}

Under these conditions it follows that%
\begin{align}
\widehat{\mathbf{\beta}}(\mathbf{y)}  & =\mathbf{\arg\max_{\beta}}\text{
}T^{-1}L_{a}(\mathbf{y(\phi}_{0}\mathbf{);\beta}%
)\overset{}{\overset{a.s.}{\rightarrow}\mathbf{\arg\max_{\beta}}\text{
}L_{\infty}(\mathbf{\phi}_{0}\mathbf{,\beta})=}\mathbf{b}(\mathbf{\phi}%
_{0}\mathbf{)}\label{betay}\\
& \text{and}\nonumber\\
\widehat{\mathbf{\beta}}(z\mathbf{)}  & =\mathbf{\arg\max_{\beta}}\text{
}T^{-1}L_{a}(\mathbf{z(\phi);\beta})\overset{}{\overset{a.s.}{\rightarrow
}\mathbf{\arg\max_{\beta}}\text{ }L_{\infty}(\mathbf{\phi,\beta})=}%
\mathbf{b}(\mathbf{\phi).}\label{betaz}%
\end{align}
Hence, in the ABC context, in which the generic parameter $\mathbf{\phi}$
represents a draw $\mathbf{\phi}^{i}$ from the prior, $p(\mathbf{\phi)}$, we
see that as $T\rightarrow\infty$ the choice criterion in (\ref{dist_mle})
approaches%
\[
d\{\mathbf{\eta}(\mathbf{y}),\mathbf{\eta}(\mathbf{z}^{i})\}=\sqrt{\left[
\mathbf{b}(\mathbf{\phi}_{0}\mathbf{)-b}(\mathbf{\phi}^{i}\mathbf{)}\right]
^{\prime}\mathbf{\Omega}\left[  \mathbf{b}(\mathbf{\phi}_{0}\mathbf{)-b}%
(\mathbf{\phi}^{i}\mathbf{)}\right]  }\leq\varepsilon.
\]
As $\varepsilon\rightarrow0$, being the relevant addition limiting condition
required in the ABC setting, we see that irrespective of the form of
$\mathbf{\Omega}$, the only values of $\mathbf{\phi}^{i}$ that will be
selected and, hence, be used to construct an estimate of the posterior
distribution, are values such that $\mathbf{b}(\mathbf{\phi}_{0}%
\mathbf{)}=\mathbf{b}(\mathbf{\phi}^{i}\mathbf{).}$ Given the assumption of
the uniqueness of the solution of $\mathbf{b(.)}$ for $\mathbf{\phi\,}$\ (or
$\mathbf{\phi}_{0}$), $\mathbf{b}(\mathbf{\phi}_{0}\mathbf{)}=\mathbf{b}%
(\mathbf{\phi}^{i}\mathbf{)}$ if and only if $\mathbf{\phi}_{0}=\mathbf{\phi
}^{i}.$ Hence, ABC produces draws that produce a degenerate distribution at
the true parameter $\mathbf{\phi}_{0}$, as required by the Bayesian
consistency property. Once again, this is despite the fact that asymptotic
sufficiency will not be achieved in the typical case in which the
approximating model is in error, a result that is analogous to the frequentist
finding of consistency for the II estimator, without full (Cramer Rao)
efficiency obtaining.

We conclude this section by citing related work in hidden Markov models (e.g.
Yildirim \textit{et al.,} 2013; Dean \textit{et al.}, 2014), in which ABC
principles that avoid summarization have been advocated. Specifically, the
difference between the observed data and the simulated pseudo-data is operated
time step by time step, as in $\prod_{t=1}^{T}\mathbb{I}_{d\{z_{t}^{i}%
,y_{t}\}\text{ }\leq\epsilon}$. This form of ABC approximation also allows for
the derivation of consistency properties (in the number of observations) of
the ABC estimates. In particular, using such a distance in the algorithm
allows for the approximation to converge to the genuine posterior when the
tolerance $\varepsilon$ goes to zero.\footnote{This is also the setting in
which Fearnhead and Prangle (2012) show that noisy ABC (Wilkinson, 2013) is
well-calibrated, i.e. converges to the relevant posterior distribution
(determined by the choice of summary statistics in this case). See also, Dean
and Singh (2011) and Dean \emph{et al.} (2014).} One problem with this
approach, however, is that the acceptance rate decreases quickly with $T$,
unless $\varepsilon$ is increasing with $T$. Jasra \textit{et al.} (2014)
provide some solutions here, but within an observation-driven model context
only. Finally, looking at the application of this form of ABC approximation in
a particle MCMC (PMCMC) setting, Jasra \textit{et al.} (2013) and Martin
\textit{et al.} (2014) establish convergence (to the exact posterior), in
connection with the alive particle filter (Le Gland and Oudjane, 2006).

\subsection{Score-based implementation}

With large computational gains, $\mathbf{\eta}(\mathbf{.})$ in (\ref{distance}%
) can be defined using the score of the auxiliary model. That is, the score
vector associated with the approximating model, when evaluated at the
simulated data, and with $\widehat{\mathbf{\beta}}(\mathbf{y)}$ substituted
for $\mathbf{\beta}$, will be closer to zero the `closer' is the simulated
data to the true. Hence, the choice criterion in (\ref{distance}) for an ABC
algorithm can be based on $\mathbf{\eta}(.)=\left.  \mathbf{S}(.;.)\right\vert
_{\mathbf{\beta=}\widehat{\mathbf{\beta}}(\mathbf{y)}}$, where
\begin{equation}
\mathbf{S}(\mathbf{z}^{i};\mathbf{\beta})=T^{-1}\frac{\partial L_{a}%
(\mathbf{z}^{i};\mathbf{\beta})}{\partial\mathbf{\beta}},\label{score}%
\end{equation}
yielding%
\begin{equation}
d\{\mathbf{\eta}(\mathbf{y}),\mathbf{\eta}(\mathbf{z}^{i})\}=\sqrt{\left[
\mathbf{S}(\mathbf{z}^{i};\widehat{\mathbf{\beta}}(\mathbf{y)})\right]
^{\prime}\mathbf{\Sigma}\left[  \mathbf{S}(\mathbf{z}^{i}%
;\widehat{\mathbf{\beta}}(\mathbf{y)})\right]  }\leq\varepsilon
,\label{dist_score}%
\end{equation}
where $\mathbf{\Sigma}$ denotes an arbitrary positive definite weighting
matrix. Implementation of ABC via (\ref{dist_score}) is faster (by many orders
of magnitude) than the approach based upon $\mathbf{\eta}%
(.)=\widehat{\mathbf{\beta}}(.)$, due to the fact that maximization of the
approximating model is required only once, in order to produce
$\widehat{\mathbf{\beta}}(.)$ from the observed data $\mathbf{y.}$ All other
calculations involve simply the \textit{evaluation}\emph{\ }of $\mathbf{S}%
(.;.)$ at the simulated data, with a numerical differentiation technique
invoked to specify $\mathbf{S}(.;.).$

Once again in line with the proof of the consistency of the relevant
frequentist (EMM) estimator, the Bayesian consistency result in Section
\ref{Theory} could be re-written in terms $\mathbf{\eta}(.)=\left.
\mathbf{S}(.;.)\right\vert _{\mathbf{\beta=}\widehat{\mathbf{\beta}%
}(\mathbf{y)}}$, upon the addition of a differentiability condition regarding
$L_{a}(\mathbf{z}^{i};\mathbf{\beta})$ and the assumption that $\mathbf{\beta
}=b(\mathbf{\phi)}$ is the unique solution to the limiting first-order
condition,\newline$\partial L_{\infty}(\mathbf{\phi,\beta})/\partial
\mathbf{\beta=}\lim_{T}T^{-1}\partial L_{a}(\mathbf{z}(\mathbf{\phi
});\mathbf{\beta})/\partial\mathbf{\beta}$ and the convergence is uniform in
$\mathbf{\beta}$. In brief, given that $\widehat{\mathbf{\beta}}%
(\mathbf{y)}\overset{a.s.}{\rightarrow}\mathbf{b}(\mathbf{\phi}_{0}\mathbf{)}%
$, as $T\rightarrow\infty$ the choice criterion in (\ref{dist_score})
approaches
\[
d\{\mathbf{\eta}(\mathbf{y}),\mathbf{\eta}(\mathbf{z}^{i})\}=\sqrt{\left[
\partial L_{\infty}(\mathbf{\phi}^{i};\mathbf{b}(\mathbf{\phi}_{0}%
\mathbf{)})/\partial\mathbf{\beta}\right]  ^{\prime}\mathbf{\Sigma}\left[
\partial L_{\infty}(\mathbf{\phi}^{i};\mathbf{b}(\mathbf{\phi}_{0}%
\mathbf{)})/\partial\mathbf{\beta}\right]  }\leq\varepsilon.
\]
As $\varepsilon\rightarrow0$, irrespective of the form of $\mathbf{\Sigma} $,
the only values of $\mathbf{\phi}^{i}$ that will be selected via ABC are
values such that $\mathbf{b}(\mathbf{\phi}_{0}\mathbf{)}=\mathbf{b}%
(\mathbf{\phi}^{i}\mathbf{)}$, which, given Assumption (A3), implies
$\mathbf{\phi}_{0}=\mathbf{\phi}^{i}.$

Hence, Bayesian consistency is maintained through the use of the score.
However, a remaining pertinent question concerns the impact on sufficiency
(or, more precisely, on \textit{the proximity to asymptotic sufficiency}) of
the use of the score instead of the MLE. In practical terms this question can
be re-phrased as: does the selection criterion based on $\mathbf{S}%
(.;.)$\ yield identical draws of $\mathbf{\phi}$ to those yielded by the
selection criterion based on $\widehat{\mathbf{\beta}}$? If the answer is yes
then, unambiguously, for large enough $T$ and for $\varepsilon\rightarrow
0$,\ the score- and MLE-based ABC criteria will yield equivalent estimates of
the exact posterior, with the accuracy of those (equivalent) estimates
dependent, of course, on the nature of the auxiliary model itself.

For any auxiliary model (satisfying identification and regularity conditions)
with unknown parameter vector $\mathbf{\beta}$, we expand the score function
in (\ref{score}), evaluated at $\widehat{\mathbf{\beta}}(\mathbf{y)}$, around
the point $\widehat{\mathbf{\beta}}(\mathbf{z}^{i}\mathbf{)}$ (with scaling
via $T^{-1}$ having been introduced at the outset in the definition of the
score in (\ref{score})),
\[
\mathbf{S}(\mathbf{z}^{i};\widehat{\mathbf{\beta}}(\mathbf{y)})=\mathbf{S}%
(\mathbf{z}^{i};\widehat{\mathbf{\beta}}(\mathbf{z}^{i}\mathbf{)}%
)+\mathbf{D}\left[  \widehat{\mathbf{\beta}}(\mathbf{y)}%
-\widehat{\mathbf{\beta}}(\mathbf{z}^{i}\mathbf{)}\right]  =\mathbf{D}\left[
\widehat{\mathbf{\beta}}(\mathbf{y)}-\widehat{\mathbf{\beta}}(\mathbf{z}%
^{i}\mathbf{)}\right]  ,
\]
where%
\begin{equation}
\mathbf{D}=T^{-1}\frac{\partial^{2}L_{a}(\mathbf{z}^{i}%
;\widetilde{\mathbf{\beta}}(\mathbf{z)})}{\partial\mathbf{\beta}%
\partial\mathbf{\beta}^{\prime}}\label{d}%
\end{equation}
and $\widetilde{\mathbf{\beta}}(\mathbf{z}^{i}\mathbf{)}$ denotes an (unknown)
intermediate value between $\widehat{\mathbf{\beta}}(\mathbf{y)}$ and
$\widehat{\mathbf{\beta}}(\mathbf{z}^{i}\mathbf{)}$. Hence, the (scaled)
criterion in (\ref{dist_score}) becomes%
\begin{align}
& \sqrt{\left[  \mathbf{S}(\mathbf{z}^{i};\widehat{\mathbf{\beta}}%
(\mathbf{y)})\right]  ^{\prime}\mathbf{\Sigma}\left[  \mathbf{S}%
(\mathbf{z}^{i};\widehat{\mathbf{\beta}}(\mathbf{y)})\right]  }\nonumber\\
& =\sqrt{\left[  \widehat{\mathbf{\beta}}(\mathbf{y)}-\widehat{\mathbf{\beta}%
}(\mathbf{z}^{i}\mathbf{)}\right]  ^{\prime}\mathbf{D}^{\prime}\mathbf{\Sigma
D}\left[  \widehat{\mathbf{\beta}}(\mathbf{y)}-\widehat{\mathbf{\beta}%
}(\mathbf{z}^{i}\mathbf{)}\right]  ^{\prime}}\leq\varepsilon.\label{new_score}%
\end{align}
Subject to standard conditions regarding the second derivatives of the
auxiliary model, the matrix $\mathbf{D}$ in (\ref{d}) will be of full rank and
as $T\rightarrow\infty$, $\mathbf{D}^{\prime}\mathbf{\Sigma D}\rightarrow$
some positive definite matrix (given the positive definiteness of
$\mathbf{\Sigma}$) that is some function of $\mathbf{\phi}^{i}.$ Hence, whilst
for any $\varepsilon>0$, the presence of $\mathbf{D}$ affects selection, as it
is a function of the drawn value $\mathbf{\phi}^{i}$ (through $\mathbf{z}^{i}%
$), as $\varepsilon\rightarrow0$, $\mathbf{\phi}^{i}$ will be selected via
(\ref{new_score}) if and only if $\widehat{\mathbf{\beta}}(\mathbf{y)}$ and
$\widehat{\mathbf{\beta}}(\mathbf{z}^{i}\mathbf{)}$ are equal. Similarly,
irrespective of the form of the (positive definite) weighting matrix in
(\ref{dist_mle}), the MLE criterion will produce these same selections. This
result pertains no matter what the dimension of $\mathbf{\beta}$ relative to
$\mathbf{\phi}$, i.e. no matter whether the true parameters are exactly or
over-identified by the parameters of the auxiliary model. This result thus
goes beyond the comparable result regarding the II/EMM estimators (see, for
e.g. Gouri\'{e}roux and Monfort, 1996), in that the equivalence is independent
of the form of weighting matrix used \textit{and} the form of identification
that prevails.

Of course in practice, ABC is implemented with $\varepsilon>0,$ at which point
the two ABC criteria will produce different draws. However, for the models
entertained in this paper, preliminary investigation has assured us that the
difference between the ABC estimates of the posteriors yielded by the
alternative criteria is negligible for small enough $\varepsilon.$\ Hence, we
proceed to operate solely with the score-based approach as the computationally
feasible method of extracting approximate asymptotic sufficiency in the state
space setting.

The actual selection and optimization of the tolerance level, $\varepsilon$,
has been the subject of intense scrutiny in the recent years (see, for
example, Marin \textit{et al}., 2011 for a detailed survey). What appears to
be the most fruitful path to the calibration of the tolerance is to firmly set
it within the realm of non-parametric statistics (Blum and Fran\c{c}ois, 2010)
as this provides proper convergence rates for the tolerance (rates that differ
between standard and noisy ABC; see Fearnhead and Prangle, 2012) and shows
that the optimal value stays away from zero for a given sample size. In
addition, the practical constraints imposed by finite computing time and the
necessity to produce an ABC sample of reasonable length lead us to follow the
recommendations found in Biau \textit{et al.} (2012), namely to analyze the
ABC approximation as a k-nearest neighbour technique and to exploit this
perspective to derive a practical value for the tolerance.

\subsection{Dimension reduction via marginal likelihood techniques}

An ABC algorithm induces two forms of approximation error. Firstly, and most
fundamentally, the use of a vector of summary statistics $\mathbf{\eta
}(\mathbf{y})$ to define the selection criterion in (\ref{distance}) means
that a simulation-based estimate of the posterior of interest is the outcome
of the exercise. Only if $\mathbf{\eta}(\mathbf{y})$ is sufficient for
$\mathbf{\phi}$ is $p(\mathbf{\phi|\eta}(\mathbf{y}))$ equivalent to the exact
posterior $p(\mathbf{\phi|y});$ otherwise the exact posterior is necessarily
estimated with error because of the analytical difference between the exact
density and the \textit{partial }posterior density $p(\mathbf{\phi|\eta
}(\mathbf{y}))$. Secondly, the partial posterior density itself,
$p(\mathbf{\phi|\eta}(\mathbf{y})),$ will be estimated with simulation error.
Critically, as highlighted by Blum (2010b), the accuracy of the
simulation-based estimate of $p(\mathbf{\phi|\eta}(\mathbf{y}))$ will be less,
all other things given, the larger the dimension of $\mathbf{\eta}%
(\mathbf{y})$. This `curse of dimensionality' obtains even when the parameter
$\mathbf{\phi}$ is a scalar, and relates solely to the dimension of
$\mathbf{\eta}(\mathbf{y}).$ As elaborated on further by Nott \textit{et al.}
(2014), this problem is exacerbated as the dimension of $\mathbf{\phi}$ itself
increases, firstly because an increase in the dimension of $\mathbf{\phi}$
brings with it a concurrent need for an increase in the dimension of
$\mathbf{\eta}(\mathbf{y})$ and, secondly, because the need to estimate a
multi-dimensional density (for $\mathbf{\phi}$) brings with it its own
problems related to dimension.

As a potential solution to the inaccuracy induced by the dimensionality of the
problem, Nott \textit{et al.} (2014) suggest allocating (via certain criteria)
a subset of the full set of summary statistics to each element of
$\mathbf{\phi}$, $\phi_{j},$ $j=1,2,....,p,$ using kernel density techniques
to estimate each marginal density, $p(\phi_{j}|\mathbf{\eta}_{j}%
(\mathbf{y})),$ and then using standard techniques to retrieve a more accurate
estimate of the joint posterior, $p(\mathbf{\phi}|\mathbf{\eta}(\mathbf{y})),$
if required. However, the remaining problem associated with the (possibly
still high) dimension of each $\mathbf{\eta}_{j}(\mathbf{y})$, in addition to
the very problem of defining an appropriate set $\mathbf{\eta}_{j}%
(\mathbf{y})$ for each $\phi_{j}$, remains unresolved. See Blum \textit{et
al.} (2013) for further elaboration on the dimensionality issue in ABC and a
review of current approaches for dealing with the problem.

The principle advocated in this paper is to exploit the information content in
the MLE of the parameters of an auxiliary model, $\mathbf{\beta}$, to yield
`approximate' asymptotic sufficiency for $\mathbf{\phi}$. Within this
framework, the dimension of $\mathbf{\beta}$ determines the dimension of
$\mathbf{\eta}(\mathbf{y})$ and the curse of dimensionality thus prevails for
high-dimensional $\mathbf{\beta.}$ However, in this case a solution is
available, at least when the dimensions of $\mathbf{\beta}$ and $\mathbf{\phi
}$ are equivalent and there is a one-to-one match between the elements of the
two parameter vectors. This is clearly so for the two cases tackled in this
paper. In the linear Gaussian model investigated as Case 1) the auxiliary
model coincides exactly with the true model, in which case $\mathbf{\beta
}=\mathbf{\phi.}$ In the stochastic volatility SSM investigated as Case 2), we
produce an auxiliary model by discretizing the latent diffusion (and
evaluating the resultant likelihood via the AUKF); hence, $\dim(\mathbf{\beta
)=}$ $\dim(\mathbf{\phi)}$ and there is a natural one-to-one mapping between
the parameters of the (true) continuous time and (approximating) discretized
models. In both of these examples then, marginalizing the auxiliary likelihood
function with respect to all parameters other than $\beta_{j}$ and then
producing the score of this function with respect to $\beta_{j}$ (as evaluated
at the marginal MLE from the observed data, $\widehat{\beta}_{j}(\mathbf{y)}%
$), yields, by construction, an obvious \textit{scalar} statistic for use in
selecting draws of $\phi_{j}$ and, hence, a method for estimating $p(\phi
_{j}\mathbf{|y}).$ If the marginal posteriors only are of interest, then all
$p$ marginals can be estimated in this way, with $p$ applications of
$(p-1)$-dimensional integration required at each step within ABC to produce
the relevant score statistics. Importantly, we do not claim here that the
`proximity' to sufficiency (for $\mathbf{\phi}$) of the vector statistic
$\mathbf{\eta(y)}$, translates into an equivalent relationship between the
score of the marginalized (auxiliary) likelihood function and the
corresponding scalar parameter, nor that the associated product density is
coherent with a joint probability distribution. If the joint posterior (of the
full vector $\mathbf{\phi}$) is of particular interest, the sort of techniques
advocated by Nott et al. (2014), amongst others, can be used to yield joint
inference from the estimated marginals.

In Section \ref{assess} we explore the benefits of marginalization, in
addition to the increase in accuracy yielded by using a score-based ABC method
(either joint or marginal) rather than an ABC algorithm based on a more
\textit{ad hoc} choice of summary statistics. However, prior to that we
provide details in the following section of the form of auxiliary model
advocated for the non-linear SSM case.

\section{The AUKF approximation for the general non-linear SSM\label{model}}

When the SSM defined in (\ref{meas_gen}) and (\ref{state_gen}) is analytically
intractable, an approximating model is needed to drive the score-based ABC
technique. In the canonical non-linear example being emphasized in the paper,
in which either $y_{t}$ or $x_{t}$ (or both) is driven by a continuous time
process, this approximation begins with the specification of a discretized
version of (\ref{meas_gen}) and (\ref{state_gen}), expressed generically using
a regression formulation as:%
\begin{align}
y_{t}  & =h_{t}\left(  x_{t},e_{t},\mathbf{\phi}\right) \label{discrete_meas}%
\\
x_{t}  & =k_{t}\left(  x_{t-1},v_{t},\mathbf{\phi}\right)
,\label{discrete_state}%
\end{align}
for $t=1,2,...,T$, where the $\left\{  e_{t}\right\}  $ and $\left\{
v_{t}\right\}  $ are assumed to be sequences of $i.i.d.$\ random variables.
This formulation is general enough to include, for example, independent random
jump components in either the measurement or state equations (subsumed under
$e_{t}$ and $v_{t}$ respectively), but does exclude cases where the nature of
the model is such that a regression formulation in discrete time is not
feasible. We note here that in producing (\ref{discrete_meas}) and
(\ref{discrete_state}) either the observation or the state variable, or both,
may need to be transformed; however, for notational simplicity we continue to
use the same symbols, $x_{t}$ and $y_{t}$, as are used to denote the variables
in the true model in (\ref{meas_gen}) and (\ref{state_gen}). The nature of the
discretization affects the functional form of $h_{t}(.)$ and $k_{t}(.).$ In
Section \ref{hest} we illustrate the method using a model in which the true
model for $y_{t}$ is already expressed in discrete time (i.e. there is no
discretization error via the measurement process) and in which a (first-order)
Euler process is initially used to approximate a square root diffusion model
for the state.

The log-likelihood function associated with the approximate model in
(\ref{discrete_meas}) and (\ref{discrete_state}) is defined by
\begin{equation}
L_{a}(\mathbf{y;\phi})=\ln p\left(  y_{1}\right)  +%
%TCIMACRO{\tsum \limits_{t=1}^{T-1}}%
%BeginExpansion
{\textstyle\sum\limits_{t=1}^{T-1}}
%EndExpansion
\ln p\left(  y_{t+1}|\mathbf{y}_{1:t}\right)  .\label{approx_like}%
\end{equation}
where $\mathbf{y}_{1:t}=(y_{1},y_{2},...,y_{t})^{\prime}.$ Only if the
approximate model is linear and Gaussian or the state variable is discrete on
a finite support, are the components used to define (\ref{approx_like})
available in closed form. Given the nature of the problem we are tackling
here, namely one in which (\ref{discrete_meas}) and (\ref{discrete_state}) are
produced as discrete approximations to a continuous time model (or, indeed,
one in which (\ref{discrete_meas}) and (\ref{discrete_state}) represent an
initial discrete-time formulation that is non-linear and/or non-Gaussian), it
is reasonable to assume that (\ref{approx_like}) cannot be evaluated exactly.
Whilst several methods (see, e.g. Simon, 2006) are available to approximate
$L_{a}(\mathbf{y;\phi})$ including, indeed, simulation-based methods such as
particle filtering and SMC, the fact that this computation is to be
\textit{embedded} within the ABC algorithm makes it essential that the
technique is both fast and numerically stable.\emph{\ }The AUKF satisfies
these criteria and, hence, is our method of choice for this illustration.

In brief, the unscented Kalman filter (UKF) is based on the theory of
unscented transformations, which is a method for calculating the moments of a
non-linear transformation of a random variable. It involves selecting a set of
points on the support of the random variable, called \textit{sigma points},
according to a predetermined and deterministic criterion. These sigma points
yield, in turn, a `cloud' of transformed points through the non-linear
function, which are then used to produce approximations to the moments of the
state variable $x_{t}$ used in implementing the KF, via simple weighted sums
of the transformed sigma points. (See Haykin, 2001, Chapter 7, for an
excellent introduction.) The \textit{augmented v}ersion of the UKF \ - the
AUKF - generalizes the filter to the case where the state and measurement
errors are non-additive, applying the principles of unscented transformations
to the \textit{augmented }state vector $s_{t}=(x_{t},\nu_{t},e_{t})^{\prime}.$
The computational burden of the filter is thus minimal - comprising the
calculation of updated sigma points for the time varying state $x_{t}$ at each
$t$, the computation of the relevant means and variances (for both $x_{t}$ and
$y_{t}$) using simple weighted sums, and the use of the usual KF up-dating
equations. Details of the specification of the sigma points, plus the steps
involved in estimating (\ref{approx_like}) are provided in the Appendix.

In implementing the score-based ABC method in a non-linear continuous time
setting, there are two aspects of the approximation used therein to consider:
1) the accuracy of the discretization; and 2) the accuracy with which the
likelihood function of the approximate (discretized) model is evaluated via
the AUKF and, hence, the accuracy of the resultant estimate of the MLE. In
addressing the first aspect, one has access to the existing literature in
which various discretized versions of continuous time models are derived, to
different orders of accuracy. With regard to the accuracy of the AUKF
evaluation of the likelihood function, we make reference to relevant results
(see, e.g. Haykin, 2001) that document the higher-order accuracy (relative,
say, to the extended KF) of the AUKF-based estimates of the mean and variance
of the filtered and predictive (for both state and observed) distributions. We
also advocate here using any transformation of the measurement (or state)
equation that renders the Gaussian approximations invoked by the AUKF more
likely to be accurate. However, beyond that, the accuracy of the Gaussian
assumption embedded in the AUKF-based likelihood estimate is
case-specific.\footnote{One potential benefit of the AUKF-based approach is
that higher-order discretizations, which produce non-linearity in the relevant
error terms could be adopted (e.g. Milstein, 1998), as the relevant sigma
points that underpin the AUKF method need only be substituted into this
additional non-linear function. We do not explore this option in the current
paper, basing the numerical demonstration in Section \ref{hest} on a
first-order Euler discretization only.}

\section{Numerical assessment of alternative ABC methods \label{assess}}

We now undertake a numerical exercise in which the accuracy of the score-based
methods of ABC (both joint and marginal) is compared with that of ABC methods
based on a set of summary statistics that are chosen without reference to a
model. We conduct this exercise firstly within the context of the linear
Gaussian model, in which case the exact score is accessible via the KF. The
results here thus provide evidence on two points of interest: 1) whether or
not accuracy can be increased by accessing the asymptotic sufficiency of the
(exact) MLE in an ABC treatment of a state space model (compared to the use of
other summary statistics); and 2) whether or not the curse of dimensionality
can be obviated via marginalization, or integration, of the exact likelihood.

In Section \ref{hest} we then assess accuracy in the typical setting in which
an approximating model is used to generate the score. The square root
volatility model is used as the example, with the approximate score produced
by evaluating the likelihood function of a discretized version of that model
using the AUKF. Whilst recognizing that the results relate to one particular
model and approximation thereof, they do, nevertheless, serve to illustrate
that score-based methods can dominate the summary statistic-based techniques
(overall), even when the approximating model used is not particularly accurate.

\subsection{Case 1: Linear Gaussian model\label{lg}}

\subsubsection{Data generation and computational details}

In this section we simulate a sample of size $T=400$ from the linear Gaussian
(LG) model in (\ref{measlg}) and (\ref{statelg}), based on the parameter
settings: $\rho=0.7$, $\delta=0.1$ and $\sigma_{v}^{2}=1$, with the
three-dimensional parameter $\mathbf{\phi}=(\rho,\delta,\sigma_{v})^{\prime}$
to be estimated. The value of the measurement error variance, $\sigma_{e}^{2}%
$, is set in order to fix the SN ratio. We compare the performance of the
(joint and marginal) score-based techniques with that of more conventional ABC
methods based on summary statistics that may be deemed to be a sensible choice
in this setting. Given the relationship between the LG model and an observable
AR(1) process, it seems sensible to propose a set of summary statistics that
are sufficient for the latter, as given in (\ref{AR1_summ_stats}). Two forms
of distances are used. Firstly, we apply the conventional Euclidean distance,
with each summary statistic also weighted by the inverse of the variance of
the values of the statistic across the ABC draws. That is, we define
\begin{equation}
d\{\mathbf{\eta}(\mathbf{y}),\mathbf{\eta}(\mathbf{z}^{i})\}=[%
%TCIMACRO{\tsum \limits_{j=1}^{5}}%
%BeginExpansion
{\textstyle\sum\limits_{j=1}^{5}}
%EndExpansion
(s_{j}^{i}-s_{j}^{obs})^{2}/var(s_{j})]^{1/2}\label{Euclid}%
\end{equation}
for ABC iteration $i=1.,2,...,R$, where $var(s_{j})$ is the variance (across
$i$) of the $s_{j}^{i}$, and $s_{j}^{obs}$ is the observed value of the $j$
statistic. Secondly, we use a distance measure proposed in Fearnhead and
Prangle (2012) which, as made explicit in Blum \textit{et al.} (2013), is a
form of dimension reduction method. We explain this briefly as follows. Given
the vector of observations $\mathbf{y}$, the set of summary statistics in
(\ref{AR1_summ_stats}) are used to produce an estimate of $E(\phi
_{j}|\mathbf{y)}$, $j=1,2,3,$ which, in turn, is used as the summary statistic
in a subsequent ABC algorithm. The steps of the procedure (as modified for
this context) described for selection of the scalar parameter $\phi_{j}$,
$j=1,2,3,$ are as follows:

\begin{enumerate}
\item Simulate $\phi_{j}^{i}$, $i=1,2,...,R$, from $p(\phi_{j})$ and,
subsequently, simulate $\mathbf{x}^{i}=(x_{1}^{i},x_{2}^{i},...,x_{T}%
^{i})^{\prime}$ from (\ref{vol}) using the exact transitions, and pseudo data,
$\mathbf{z}^{i}$ using the conditional Gaussian form of $p(\mathbf{z}%
|\mathbf{x}).$

\item For $\mathbf{z}^{i}$, $i=1,2,...,R$, calculate
\begin{equation}
\mathbf{s}^{i}=\left[  s_{1}^{i},s_{2}^{i},s_{3}^{i},s_{4}^{i},s_{5}%
^{i}\right]  ^{\prime}\label{fp_stats}%
\end{equation}

\item Define $\mathbf{\phi}_{j}=(\phi_{j}^{1},\phi_{j}^{2},...,\phi_{j}%
^{R})^{\prime},$ $\mathbf{X}=\left[
\begin{array}
[c]{cccc}%
1 & 1 & \cdots & 1\\
\mathbf{s}^{1} & \mathbf{s}^{2} & \cdots & \mathbf{s}^{R}%
\end{array}
\right]  ^{\prime}$ and
\[
\mathbf{\phi}_{j}=E[\mathbf{\phi}_{j}|\mathbf{Z}]+\mathbf{e}=\mathbf{X}\left[
\begin{array}
[c]{c}%
\alpha\\
\mathbf{\beta}%
\end{array}
\right]  +\mathbf{e,}%
\]
where $\mathbf{Z}=\left[  \mathbf{z}^{1},\mathbf{z}^{2},...,\mathbf{z}%
^{T}\right]  $ and $\mathbf{\beta}$ is of dimension $(5\times1).$

\item Use OLS to estimate $E[\mathbf{\phi}_{j}|\mathbf{Z}]$\textbf{\ }as
$\widehat{E}[\mathbf{\phi}_{j}|\mathbf{Z}]=\widehat{\alpha}+\left[
\begin{array}
[c]{cccc}%
\mathbf{s}^{1} & \mathbf{s}^{2} & \cdots & \mathbf{s}^{R}%
\end{array}
\right]  ^{\prime}\widehat{\mathbf{\beta}}$

\item Define:%
\[
\eta(\mathbf{z}^{i})=\widehat{E}(\phi_{j}|\mathbf{z}^{i})=\widehat{\alpha
}+\mathbf{s}^{i^{\prime}}\widehat{\mathbf{\beta}}\text{ and }\eta
(\mathbf{y})=\widehat{E}(\phi_{j}|\mathbf{y})=\widehat{\alpha}+\mathbf{s}%
^{obs^{\prime}}\widehat{\mathbf{\beta}},
\]
where $\mathbf{s}^{obs}$ denotes the vector of summary statistics in
(\ref{fp_stats}) calculated from the vector of observed returns, and use:%
\begin{equation}
d\{\eta(\mathbf{y}),\eta(\mathbf{z}^{i})\}=\left\vert \widehat{E}(\phi
_{j}|\mathbf{y})-\widehat{E}(\phi_{j}|\mathbf{z}^{i})\right\vert =\left\vert
\mathbf{s}^{i^{\prime}}\widehat{\mathbf{\beta}}-\mathbf{s}^{obs^{\prime}%
}\widehat{\mathbf{\beta}}\right\vert \label{fp}%
\end{equation}
as the selection criterion for $\phi_{j}$ at each iteration $i$.
\end{enumerate}

The joint score-based method uses the distance measure in (\ref{dist_score}),
but with the score in this case computed from the \textit{exact }model,
evaluated using the KF. The weighting matrix $\mathbf{\Sigma}$ is set equal to
the Hessian-based estimate of the covariance matrix of the (joint)\ MLE
estimator of $\mathbf{\beta=\phi}$, evaluated at the MLE computed from the
observed data, $\widehat{\mathbf{\phi}}(\mathbf{y).}$ The marginal score-based
method for estimating the marginal posterior for the $jth$ element of
$\mathbf{\phi}$, $\phi_{j}$, $j=1,2,3,$ is based on the distance%
\begin{equation}
d\{\mathbf{\eta}(\mathbf{y}),\mathbf{\eta}(\mathbf{z}^{i})\}=\left\vert
S(\mathbf{y};\phi_{j}\mathbf{)}-S(\mathbf{z}^{i};\phi_{j})\right\vert
,\label{dist_marg}%
\end{equation}
where%
\begin{equation}
S(\mathbf{z}^{i};\phi_{j})=T^{-1}\frac{\partial L(\mathbf{z}^{i};\phi_{j}%
)}{\partial\phi_{j}}\label{marg_score}%
\end{equation}
and $L(\mathbf{z}^{i};\phi_{j})$ is produced by integrating (numerically) the
exact likelihood function (evaluated via the KF) with respect to all
parameters other than $\phi_{j}$, and taking the logarithm.

\subsubsection{Numerical results for the LG model\label{lg results}}

\begin{figure}[ptb]%
\centering
\caption{Marginal posterior densities for the three parameters in the linear
Gaussian state space model in (\ref{measlg}) and (\ref{statelg}). As per the
key, the graphs reproduced are the exact posterior, in addition to the four
ABC-based estimates, as detailed in the text.\bigskip}%
\includegraphics[
natheight=10.824900in,
natwidth=14.222700in,
height=5.4405in,
width=6.8554in
]%
{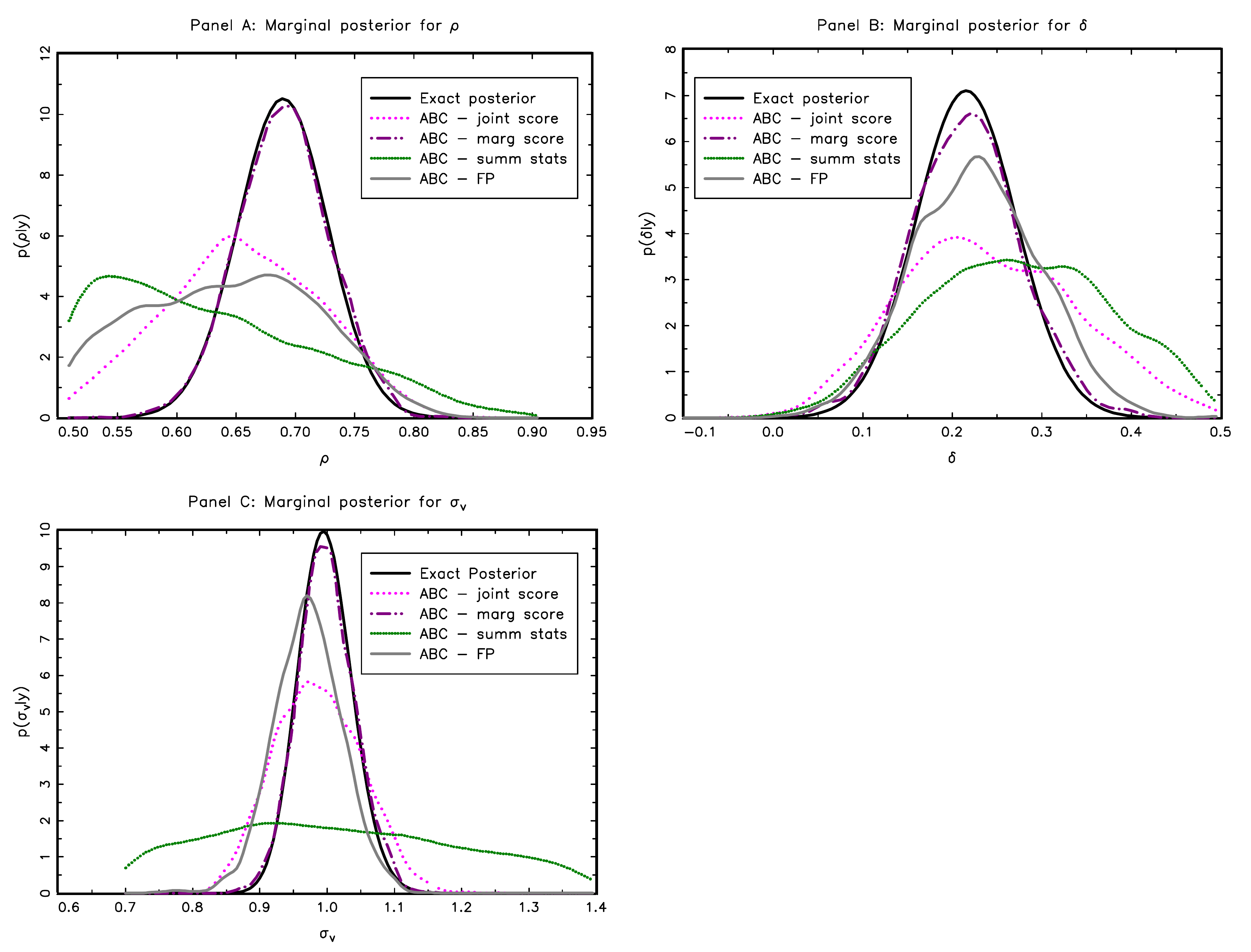}%
\label{fig_lg}%
\end{figure}

We produce results that compare the performance of the four different methods:
the joint score-based ABC (`ABC-joint score' in all figures); the marginal
score-based ABC (`ABC-marg score'); the summary statistic-based ABC using the
Euclidean metric in (\ref{Euclid}) (`ABC-summ stats'); and the approach of
Fearnhead and Prangle (2012) based on the metric in (\ref{fp}) (`ABC-FP').
Marginal density estimates are produced initially for a single run of ABC,
based on 50,000 replications of the accept/reject algorithm detailed in
Section 2.1, and with $\varepsilon$ defined as the 5th percentile of the
50,000 draws. The true data is generated from a process in which the SN ratio
is high (i.e. $\left[  \sigma_{v}^{2}/(1-\rho^{2})\right]  /\sigma_{e}^{2}%
=20$), with the true marginal posteriors computed by normalizing the
likelihood function evaluated using the KF (and multiplied by a uniform
prior), then marginalizing using deterministic integration over a very fine
grid for $\mathbf{\phi.}$ The score-based methods also use the exact
likelihood function to compute the score. All three sets of marginal
posteriors, for $\rho$, $\delta$ and $\sigma_{v}$, are produced in Figure
\ref{fig_lg}, Panels A, B and C respectively. We then summarize the results
for 100 replications of ABC, using box plots, in Figures 2 to 4. We produce
estimates of the $5th$, $25th$, $50th$, $75th$ and $95th$ percentiles of each
true (KF-based) posterior density, with the exact percentile represented by
the horizontal dotted line. The short-hand notation used to denote the form of
ABC method corresponds to that used in Figure 1.%

%EndExpansion

\begin{figure}[ptb]%
\centering
\caption{Box plots for 100 replications of ABC, with each replication based on 50,000 draws. Estimation of the marginal posterior density of the state intercept parameter $\rho$ in the linear Gaussian model in (\ref{measlg}) and (\ref{statelg}). $T=400.$ High SN ratio. The percentiles for the exact marginal posterior are represented by the horizontal dotted lines; and the four ABC methods referenced correspond to those described in the text.}
\includegraphics[
natheight=6.280300in,
natwidth=12.800100in,
height=3.1678in,
width=7.0681in
]%
{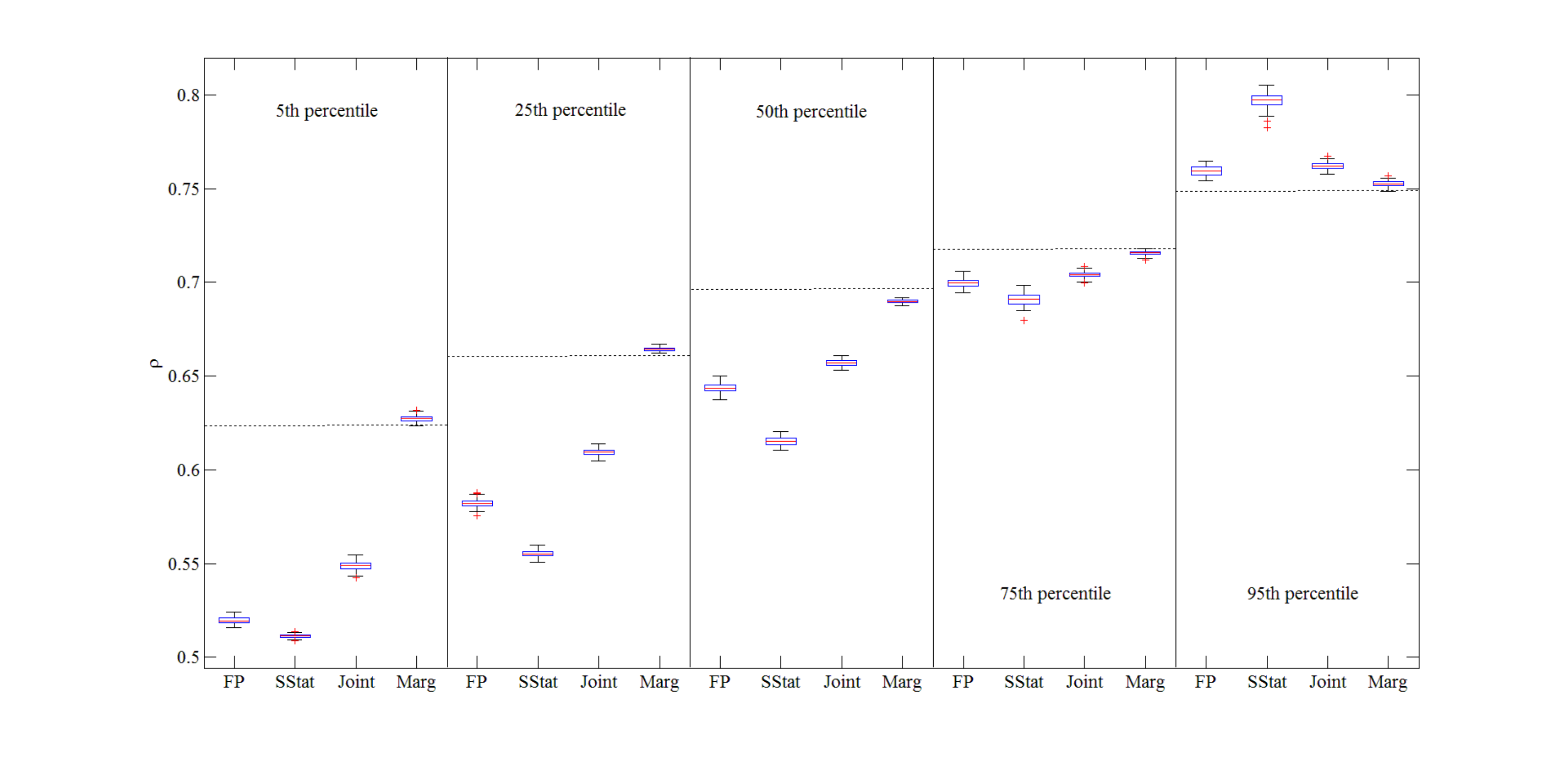}%
\label{bp_rho_high}%
\end{figure}

\begin{figure}[ptb]%
\centering
\caption{Box plots for 100 replications of ABC, with each replication based on 50,000 draws. Estimation of the marginal posterior density of the state intercept parameter $\delta$ in the linear Gaussian model in (\ref{measlg}) and (\ref{statelg}). $T=400.$ High SN ratio. The percentiles for the exact marginal posterior are represented by the horizontal dotted lines; and the four ABC methods referenced correspond to those described in the text.}
\includegraphics[
natheight=6.280300in,
natwidth=12.800100in,
height=3.1687in,
width=7.0681in
]%
{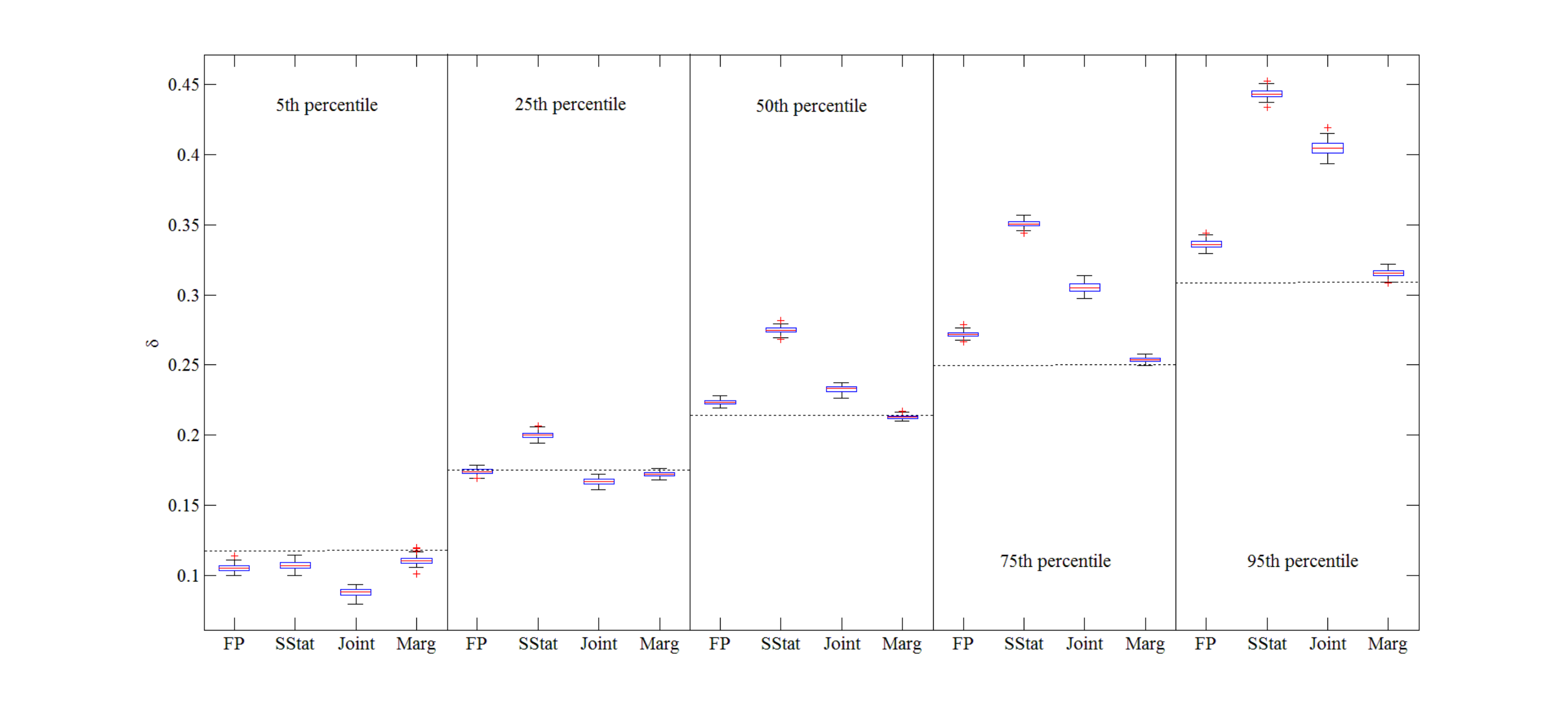}%
\label{bp_d_high}%
\end{figure}

\begin{figure}[ptb]%
\centering
\caption{Box plots for 100 replications of ABC, with each replication based on 50,000 draws. Estimation of the marginal posterior density of the state intercept parameter $\sigma_{v}$ in the linear Gaussian model in (\ref{measlg}) and (\ref{statelg}). $T=400.$ High SN ratio. The percentiles for the exact marginal posterior are represented by the horizontal dotted lines; and the four ABC methods
referenced correspond to those described in the text.}
\includegraphics[
natheight=6.280300in,
natwidth=12.800100in,
height=3.1687in,
width=7.0681in
]%
{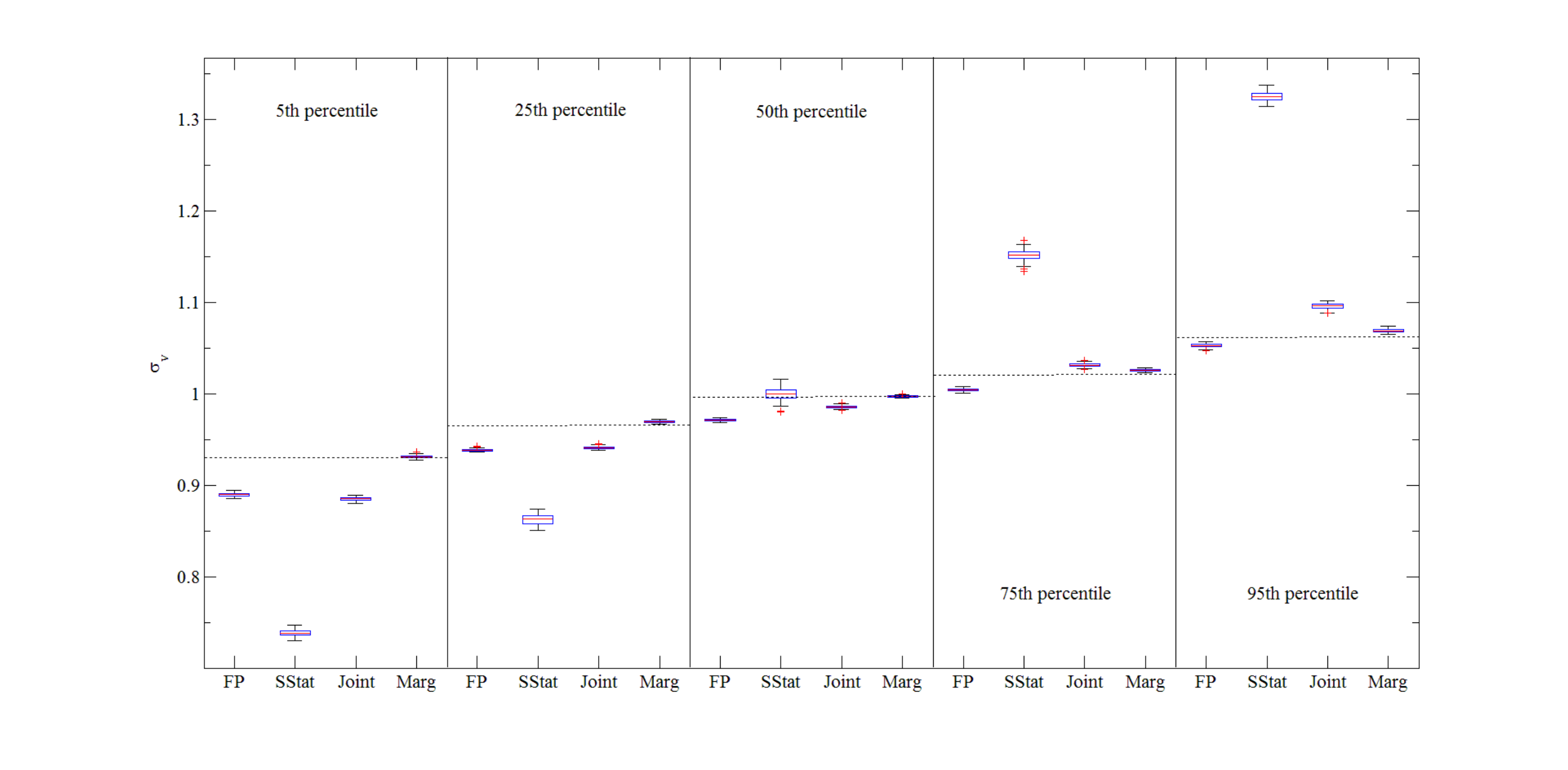}%
\label{bp_sig_high}%
\end{figure}

Figure \ref{fig_lg}, Panels A to C, highlight graphically \ - for a single ABC
run - the performance of the four ABC methods in reproducing the true marginal
posteriors. The two most notable features of all three plots are: i) the
remarkable accuracy of the marginal score approach; and ii) the very poor
performance of the summary statistic approach based on the Euclidean metric.
These features are replicated across the 100 runs of ABC, as evidenced by
Figures \ref{bp_rho_high} to \ref{bp_sig_high}. For all three parameters, the
marginal score approach produces percentile estimates that are extremely
accurate in terms of both mean location and spread. For the parameter $\rho$
the joint score approach is the next most accurate technique, followed by the
FP method. For the parameter $\sigma_{v}$ there is little to choose between
the latter two methods, both of which produce quite accurate estimates of the
true percentiles, although both being still less accurate than the marginal
score technique. For $\delta$ the FP method tends to outperform the joint
score approach, but with neither method competing with the accuracy yielded by
the marginal score. For all three parameters, the summary statistic approach
based on the Euclidean metric produces the poorest results overall, despite
the high SN ratio.

Further results (not documented here, for reasons of space) demonstrate that
decreasing the SN ratio has no qualitative effective on the performance of the
score-based results, including the marked accuracy of the marginal score
method. This robustness of the score methods to the SN ratio is to be
anticipated, give that the likelihood function for the full state space model
has been used to generate the matching statistics. The impact of the change in
the SN ratio on the summary statistic-based methods is not uniform, with there
certainly being no clear tendency for the results the worsen. This suggests
that the accuracy with which $p(\phi_{j}\mathbf{|\eta}(\mathbf{y}))$ itself is
estimated (and, hence, the issue of dimensionality), rather than the
relationship between $p(\phi_{j}\mathbf{|\eta}(\mathbf{y}))$ and $p(\phi
_{j}\mathbf{|y})$ (and, hence the `closeness' of $\mathbf{\eta}(\mathbf{y})$
to sufficiency), remains a dominant influence on final accuracy, no matter
what the SN value.

\subsection{Case 2: Heston stochastic volatility (SV) model \label{hest}}

\subsubsection{Data generation and computational details}

Given the critical role played by volatility in asset pricing, portfolio
management and the calculation of risk measures, a large segment of the
empirical finance literature has been devoted to the construction and analysis
of volatility models. Three decades of empirical studies have demonstrated
that the \textit{constant} volatility feature of a geometric Brownian motion
process for an asset price is inconsistent with both the observed
time-variation in return volatility and the non-Gaussian characteristics of
empirical distributions of returns; see Bollerslev, Chou and Kroner (1992) for
a review. Empirical regularities documented in the option pricing literature,
most notably implied volatility `smiles', are also viewed as evidence that
asset prices deviate from the geometric Brownian motion assumption underlying
the Black and Scholes (1973) option price; see, for example, Bakshi \textit{et
al.} (1997) and Lim \textit{et al.} (2005), or Garcia \emph{et al. }(2010) for
a recent review.

In response to these now well-established empirical findings, many alternative
time-varying volatility models have been proposed, with continuous time
stochastic volatility (SV) models - often augmented by random jump processes -
being particularly prominent of late. This focus on the latter form of models
is due, in part, to the availability of (semi-) closed-form option prices,
with variants of the SV model of Heston (1993) becoming the workhorse of the
empirical option pricing literature (e.g. Eraker, 2004, Forbes \textit{et al.}
2007, Broadie \textit{et al.} 2007, Johannes \textit{et al.} 2009), and MCMC
and particle filtering techniques the typical numerical methods of choice. It
is of interest, therefore, to assess the performance of the proposed ABC
method in the context of this form of model.

We adopt here the simplest version of Heston `square root' SV model, given by:%
\begin{align}
r_{t}  & =\sqrt{V_{t}}\epsilon_{t}\label{return}\\
dV_{t}  & =\left(  \delta-\alpha V_{t}\right)  dt+\sigma_{v}\sqrt{V_{t}}%
dW_{t},\label{vol}%
\end{align}
where $r_{t}$ denotes the (demeaned) logarithmic return on an asset price over
period (day say) $t$, $\epsilon_{t}\sim i.i.d.N(0,1)$, $W_{t}$ is a standard
Wiener process, and the restriction $2\delta\geq\sigma_{v}^{2}$ ensures the
positivity of the stochastic variance $V_{t}$ $(=x_{t}$ in our previous
generic notation). For $\alpha>0$, \thinspace$V_{t}$ is mean reverting and as
$t\rightarrow\infty$ the variance approaches a steady state gamma
distribution, with $E[V_{t}]=\delta/\alpha$ and $var(V_{t})=\sigma_{v}%
^{2}\delta/2\alpha^{2}.$ The transition density for $V_{t}$, conditional on
$V_{t-1}$, is%
\begin{equation}
p(V_{t}|V_{t-1})=c\exp(-u-v)\left(  \frac{v}{u}\right)  ^{q/2}I_{q}%
(2(uv)^{1/2}),\label{non-central}%
\end{equation}
where $c=2\alpha/\sigma_{v}^{2}(1-\exp(-\alpha))$, $u=cV_{t-1}\exp(-\alpha)$,
$v=cV_{t}$, $q=\frac{2\delta}{\sigma_{v}^{2}}-1$, and $I_{q}(.)$ is the
modified Bessel function of the first kind of order $q.$ The conditional
distribution function is non-Central chi-square, $\chi^{2}(2cV_{t};2q+2,2u)$,
with $2q+2$ degrees of freedom and non-centrality parameter $2u$. The
discrete-time model for $r_{t}$ in (\ref{return}) can be viewed as a
discretized version of a diffusion process for returns (or returns can be
viewed as being \textit{inherently} discretely observed) whilst we retain the
diffusion model for the latent variance. That is, we eschew the discretization
of the variance process that would typically be used, with simulation of
$V_{t}$ in Step 2 of the ABC algorithm occurring exactly, through the
treatment of the $\chi^{2}(2cV_{t};2q+2,2u)$ random variable as a composition
of central $\chi^{2}(2q+2+2j)$ and $j\sim Poisson(u)$ variables.\footnote{As
noted earlier, even in the typical case in which the transition distribution
of the diffusion is unknown, an arbitrarily fine discretization, limited only
by computer power, effectively enables draws from the exact diffusion to be
produced.}

For the purpose of this illustration we set the parameters in (\ref{vol}) to
values that produce simulated values of both $r_{t}$ and $V_{t}$ that match
the characteristics of (respectively) daily returns and daily values of
realized volatility (constructed from 5 minute returns) for the S\&P500 stock
index over the 2003-2004 period, namely%
\begin{equation}
\rho=1-\alpha=0.92;\text{ }\delta=0.0024;\text{ }\sigma_{v}%
=0.062.\label{param_values}%
\end{equation}
This relatively calm period in the stock market is deliberately chosen as a
reference point, as the inclusion of price and volatility jumps, and/or a
non-Gaussian conditional distribution in the model would be an empirical
necessity for any more volatile period, such as that witnessed during the
recent 2008/2009 financial crisis.

In order to implement the score-based ABC method, using the AUKF algorithm to
evaluate an auxiliary likelihood, we invoke the following discretization,
based on (an exact) transformation of the measurement equation and an Euler
approximation of the state equation,%
\begin{align}
\ln(r_{t}^{2})  & =\ln(V_{t})+\ln(\varepsilon_{t}^{2})\nonumber\\
& \Rightarrow\nonumber\\
y_{t}  & =\ln(V_{t})+e_{t}\label{trans}\\
& \nonumber\\
V_{t}  & =\delta+\rho V_{t-1}+\sigma_{v}\sqrt{V_{t-1}}v_{t},\label{v_state}%
\end{align}
where $v_{t}$ is treated as a truncated Gaussian variable with lower bound,%
\begin{equation}
v_{t}>\frac{-(\delta+\rho V_{t-1})}{\sigma_{v}\sqrt{V_{t-1}}}.\label{v_trunc}%
\end{equation}
Directing the reader to the Appendix for the detailed outline of the AUKF
approach, we note that sigma points that span the support of $e_{t}$ are
defined by calculating $E(e_{t})$ and $var(e_{t})$ using deterministic
integration and the closed form of $p(e_{t}),$ with the specification of
$a_{e}=b_{e}=\sqrt{3}$ adopted for convenience. Those for $V_{t}$,
$t=0,1,....,T,$ are defined as:%
\[
V_{t}^{1}=E(V_{t}^{1}|.);\text{ }V_{t}^{2}=E(V_{t}^{1}|.)+\sqrt{3}%
\sqrt{var(V_{t}^{1}|.)};\text{ }V_{t}^{3}=0.00001,
\]
where $E(V_{t}^{1}|.)$ and $var(V_{t}^{1}|.)$ respectively denote the mean and
variance of the relevant distribution of $V_{t}$ (marginal, filtered or
predictive, depending on the particular step in the AUKF algorithm). The sigma
points for (\ref{v_trunc}) are then defined using the mean and variance of the
truncated normal distribution, with the value of $V_{t-1}$ in (\ref{v_trunc})
represented using the relevant sigma point for $V_{t-1}$, and (with reference
to the Appendix) $a_{v}=b_{v}=\sqrt{3}$ specified.

In order to evaluate the accuracy of the estimate of the posterior produced
using the ABC method, we produce the exact joint posterior distribution for
$\mathbf{\phi}=(\rho,\delta,\sigma_{v})^{\prime}$ via the deterministic
non-linear filtering method of Ng \textit{et al.} (2013). In brief, this
method represents the recursive filtering and prediction distributions used to
define the likelihood function as the numerical solutions of integrals defined
over the support of $e_{t}$ in (\ref{trans}), with deterministic integration
used to evaluate the relevant integrals, and the \textit{exact} transitions in
(\ref{non-central}) used in the specification of the filtering and up-dating
steps.\footnote{A numerically efficient and stable algorithm for evaluating
the transitions densities for the exact model - in which a non-central
chi-squared density is represented as an infinite mixture of Poisson and
central chi-squared densities - is used. This enables the numerical problems
typically associated with the direct computation of Bessel functions to be
avoided.} Whilst lacking the general applicability of the ABC-based method
proposed here, this deterministic filtering method is ideal for the particular
model used in this illustration, and can be viewed as producing a very
accurate estimate of the exact density, without any of the simulation error
that would be associated with an MCMC-based comparator, for instance. We refer
the reader to Ng \textit{et al.} for more details of the technique; see also
Kitagawa (1987).\footnote{We note that the application of this filter in Ng
\textit{et al. }is to a non-parametric representation of $e_{t}.$ In the
current setting, in which $e_{t}$ is specifed parametrically, the known form
of the distribution of $e_{t}$ is used directly in the evaluation of the
relevant integrals. We refer the reader to Section 2.2. of that paper for a
full description of the algorithm. Preliminary experimentation with the number
of grid points used in the deterministic integration was undertaken in order
to ensure that the resulting estimate of the likelihood function/posterior
stabilized, with 100 grid points underlying the final results documented
here.} The likelihood function, evaluated via this method, is then multiplied
by a uniform prior that imposes the restrictions: $0<\rho<1,$ $\delta>0$,
$\sigma_{v}^{2}>0$ and $2\delta\geq\sigma_{v}^{2}$. The three marginal
posteriors are then produced via deterministic numerical integration (over the
parameter space), with a very fine grid on $\mathbf{\phi}$ being used to
ensure accuracy.

We compare the auxiliary model-based ABC technique with ABC approaches based
on summary statistics, as discussed in the linear Gaussian context in Section
\ref{lg}. For want of a better choice we use the vector of statistics given
attention therein, as well as the two alternative distance measures described
there. We also compute marginal posterior densities estimated using the
AUKF-based approximation (of the likelihood) itself. That is, using the AUKF
to evaluate the likelihood function associated with the discretized model in
(\ref{trans}) and (\ref{v_state}) and normalizing (using a uniform prior)
produces an approximation of the posterior which can, in principle, be invoked
as an approximation in its own right, independently of its subsequent use as a
score generator with an ABC algorithm. Finally, we compute the marginal
posteriors (again, based on a uniform prior) using the likelihood function of
the Euler approximation in (\ref{trans}) and (\ref{v_state}) evaluated using
the Ng \textit{et al.} (2013) filtering method, with Gaussian transitions used
in the filtering and up-dating steps. When normalized, this density can be
viewed as the quantity that a typical MCMC scheme (as based on the equivalent
prior) would be targeting, given that the tractability of Gaussian
approximations to the transitions would typically be exploited in structuring
an MCMC algorithm.

As a concluding note on computational matters, we re-iterate that the time
taken to evaluate the AUKF-based approximate likelihood function at any point
in the parameter space is roughly comparable to that required for KF
evaluation, thus rendering it a feasible method to be inserted within the ABC
algorithm. In contrast, evaluation of the Euler-based likelihood via the Ng
\textit{et al. }(2013) technique, whilst producing, in the main, (as will be
seen in the following section) a more accurate estimate of the exact posterior
than the AUKF method, is many orders of magnitude slower and, hence, simply
infeasible as a score generator within ABC.

\subsubsection{Numerical results for the Heston SV model}

In order to abstract initially from the impact of dimensionality on the ABC
methods, we first report results for each single parameter of the Heston
model, keeping the remaining two parameters fixed at their true values. Three
ABC-based estimates of the relevant exact (univariate) posterior, invoking a
uniform prior, are produced in this instance. Three matching statistics are
used, respectively: 1) the (uni-dimensional) auxiliary score based on the
approximating model (ABC-score); 2) the summary statistics in
(\ref{AR1_summ_stats}), matched via the Euclidean distance measure in
(\ref{Euclid}) (ABC-summ stats); and 3) the summary statistics in
(\ref{AR1_summ_stats}), matched via the FP distance measure in (\ref{fp})
(ABC-FP). We produce representative posterior (estimates) in each case, to
give some visual idea of the accuracy (or otherwise) that is achievable via
the ABC methods. We then summarize accuracy by reporting the average (over the
100 runs) of the root mean squared error (RMSE) of each ABC-based estimate of
the exact posterior for a given parameter, computed as:%
\begin{equation}
RMSE=\sqrt{\frac{1}{G}%
%TCIMACRO{\tsum \limits_{g=1}^{G}}%
%BeginExpansion
{\textstyle\sum\limits_{g=1}^{G}}
%EndExpansion
(\widehat{p}_{g}-p_{g})^{2}},\label{rmse_1}%
\end{equation}
where $\widehat{p}_{g}$ is the ordinate of the ABC density estimate and
$p_{g}$ the ordinate of the exact posterior density, at the $gth$ grid-point
used to produce the plots.\footnote{The associated box plots are also
available on request, but were not included here due to space considerations.}
All single parameter results are documented in Panel A of Table 1. We also
tabulate there, as benchmarks of a sort, the RMSEs associated with the
(one-off) AUKF- and Euler-based approximations of each univariate density.%

%TCIMACRO{\FRAME{fhFO}{6.8554in}{4.8983in}{0pt}{\Qct{Posterior densities for
%each single unknown parameter of the model in (\ref{trans}) and (\ref{v_state}%
%), with the other two parameters set to their true values. As per the key, the
%graphs reproduced are the exact, Euler- and AUKF-based approximations, in
%addition to the three ABC-based estimates. Both the exact and Euler posteriors
%are evaluated using the grid-based non-linear filter of Ng \QTR{em}{et al.
%}(2013).}}{\Qlb{SQ_fig}}{sq_uni_paper_graph.eps}%
%{\special{ language "Scientific Word";  type "GRAPHIC";  display "USEDEF";
%valid_file "F";  width 6.8554in;  height 4.8983in;  depth 0pt;
%original-width 14.2227in;  original-height 10.8249in;  cropleft "0";
%filename 'SQ_uni_paper_graph.eps';file-properties "XNPEU";}} }%
%BeginExpansion
\begin{figure}[h]%
\centering
\caption{Posterior densities for each single unknown parameter of the model in
(\ref{trans}) and (\ref{v_state}), with the other two parameters set to their
true values. As per the key, the graphs reproduced are the exact, Euler- and
AUKF-based approximations, in addition to the three ABC-based estimates. Both
the exact and Euler posteriors are evaluated using the grid-based non-linear
filter of Ng \emph{et al. }(2013).}%
\includegraphics[
natheight=10.824900in,
natwidth=14.222700in,
height=4.8983in,
width=6.8554in
]%
{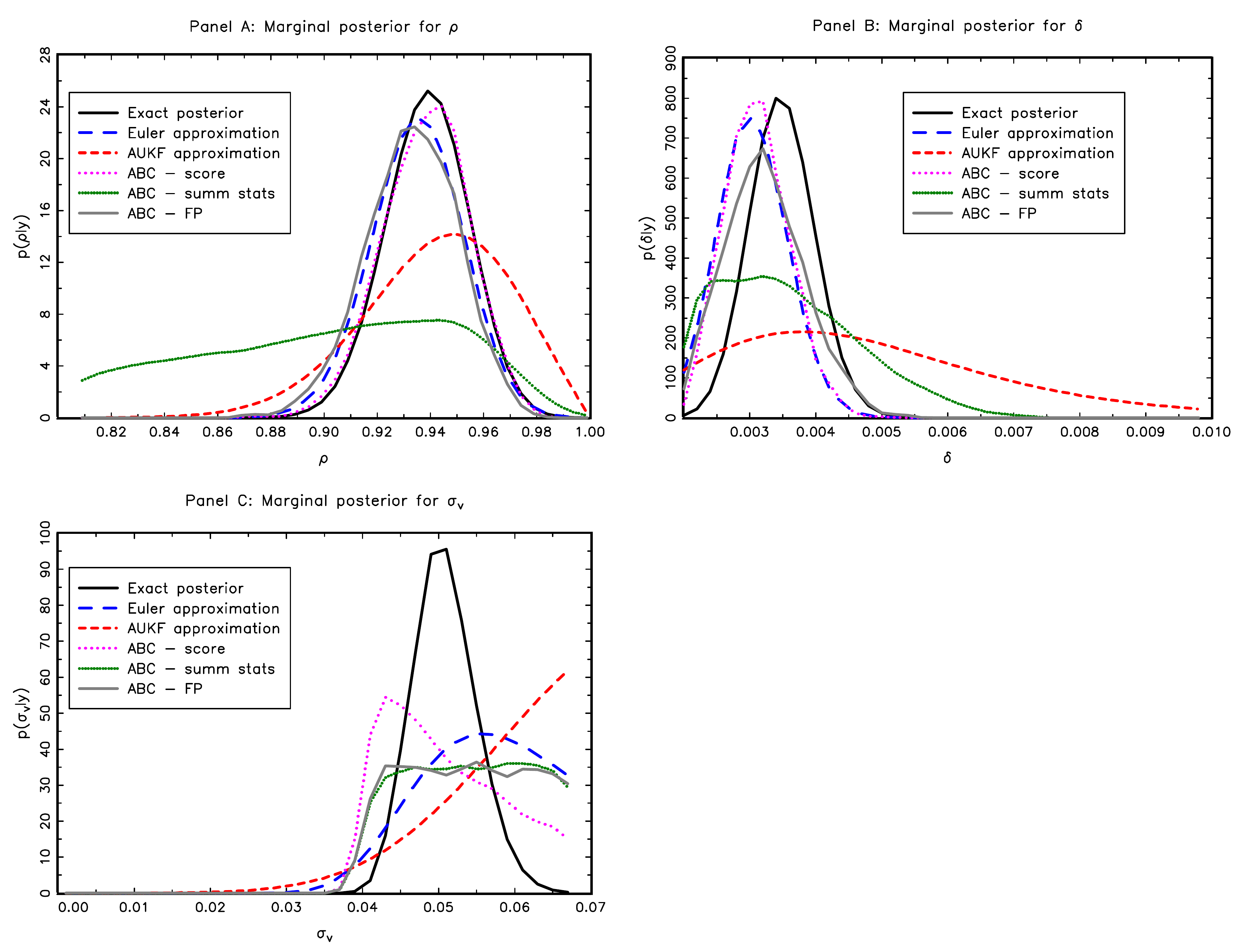}%
\label{SQ_fig}%
\end{figure}
%EndExpansion

Figure \ref{SQ_fig}, Panel A reproduces the exact posterior of (the single
unknown parameter) $\rho$, the posteriors associated with the AUKF- and
Euler-based approximations, and the three ABC-based estimates. As is clear,
the AUKF-based approximation is reasonably inaccurate, in terms of replicating
the location and shape of the exact posterior - an observation that is
interesting in its own right, given the potential for such a simple and
computationally efficient approximation method to be used to evaluate
likelihood functions (and posterior distributions) in non-linear state space
models such as the one under consideration. However, once the approximation is
embedded within an ABC scheme, in the manner described in Section \ref{model},
the situation is altogether different, with the (pink) dotted line (denoted by
`ABC-score' in the key) providing a remarkably accurate estimate of the exact
posterior, using only 50,000 replications of the simplest rejection-based ABC
algorithm, and fifteen minutes of computing time on a desktop computer. It is
worth noting that the ABC-based estimate is also more accurate in this case
than the Euler approximation, where we highlight, once again, that production
of the latter still requires the application of the much more computationally
burdensome non-linear filtering method. Most notably, the ABC method based on
the summary statistics, combined using a Euclidean distance measure, performs
very badly, although the dimensional reduction technique of Fearnhead and
Prangle (2012), applied to this same set of summary statistics, yields a
reasonable estimate of the exact posterior in this instance.

Comparable graphs are produced for the single parameters $\delta$ and
$\sigma_{v}$ in Panels B and C respectively of Figure \ref{SQ_fig}, with the
remaining pairs of parameters ($\rho$ and $\sigma_{v}$, and $\rho$ and
$\delta$ respectively) held fixed at their true values. In the case of
$\delta$, the score-based method arguably provides the best representation of
the shape of the exact posterior, despite being slightly inaccurate in terms
of location. The Fearnhead and Prangle (2012) method also provides a
reasonable estimate, whilst the summary statistic approach using the Euclidean
distance, once again performs very poorly. For the parameter $\sigma_{v}$,
\textit{only }the score based method yields a density with a well-defined
shape, with the two summary statistic-based techniques essentially producing
uniform densities that reflect little more than the restricted support imposed
on the parameter draws. Interestingly, the Euler approximation itself provides
a quite poor representation of the exact marginal, a result which has not, as
far as we know, been remarked upon in the literature, given that a typical
MCMC scheme (as noted above) would in fact be targeting the Euler density
itself, as the best representation of the true model - the exact posterior
remaining unaccessed due to the difficulty of devising an effective MCMC
scheme which uses the exact transitions. For neither $\delta$ nor $\sigma_{v}$
is the AUKF approximation \textit{itself} particularly accurate, despite the
fact that it respects the non-linearity in the true state space model.

The RMSE results recorded in Panel A of Table 1 confirm the qualitative nature
of the single-run graphical results. For $\delta$ and $\sigma_{v}$, all three
ABC-based estimates are seen to produce lower RMSE values (sometimes an order
of magnitude lower) than the AUKF approximation, indicating that \textit{any}
of the ABC procedures would yield gains over the use of the unscented
filtering method itself. For $\rho$, the AUKF approximation is better than
that of the summary statistic-based ABC estimate, but in part because the
latter is so poor. For $\rho$ and $\sigma_{v}$, the score-based ABC method is
the most accurate, and is most notably \textit{very} precise for the case of
the persistence parameter $\rho$. For the parameter $\delta$ the FP method is
the most accurate according to this measure although, as indicated by the
nature of the graphs in Panel B of Figure \ref{SQ_fig}, this result does tend
to understate the ability of the score-based method to capture the basic shape
of the exact posterior. For $\rho$ and $\delta$, the (Euclidean)
summary-statistic method is an order of magnitude more inaccurate than the
other two ABC methods, whilst also exhibiting no ability to identify the shape
of the true posterior in the case of $\sigma_{v}.$ Once again as is consistent
with the graphs in Figure \ref{SQ_fig}, the Euler approximation for $\rho$ is
reasonably accurate, but not as accurate as the score-based ABC estimate. For
$\delta$ and $\sigma_{v},$ the Euler approximation, whilst being more accurate
that the AUKF approximation, is dominated by all three ABC estimates.%

%TCIMACRO{\TeXButton{begin_landscape}{\begin{landscape}}}%
%BeginExpansion
\begin{landscape}%
%EndExpansion
%

%TCIMACRO{\TeXButton{B}{\begin{table}[tbp] \centering}}%
%BeginExpansion
\begin{table}[tbp] \centering
%EndExpansion
\label{Table1}

Table 1: RMSE of an estimated marginal and the exact marginal: average RMSE
value over multiple runs of ABC using 50,000 replications. `Score' refers to
the ABC method based on the score of the AUKF model; `SS' refers to the ABC
method based on a Euclidean distance for the summary statistics in
(\ref{AR1_summ_stats}); `FP' refers to the Fearnhead and Prangle ABC method,
based on the summary statistics in (\ref{AR1_summ_stats}). For the single
parameter case, the (single) score method is documented in the row denoted by
`ABC-Marginal Score', whilst in the multi-parameter case, there are results
for both the joint and marginal score methods. The RMSE of the AUFK and Euler
approximations (computed once only, using the observed data) are recorded as
benchmarks, in the top two rows of each panel. For the single and dual
parameter cases, 100 runs of ABC were used to produce the results, whilst for
the three parameter case, 50 runs were used. The smallest RMSE figure in each
column is highlighted in bold.%

\begin{tabular}
[c]{lllc|l|l|ll|lll}\hline\hline
& \multicolumn{3}{c|}{Panel A} & \multicolumn{2}{|c|}{Panel B} &
\multicolumn{2}{|c|}{Panel C} & \multicolumn{3}{|c}{Panel D}\\
& \multicolumn{3}{c|}{One unknown} & \multicolumn{2}{|c|}{Two unknowns} &
\multicolumn{2}{|c|}{Two unknowns} & \multicolumn{3}{|c}{Three unknowns}%
\\\cline{2-11}
& $\rho$ & $\delta$ & $\sigma_{v}$ & $\rho$ & $\sigma_{v}$ & $\rho$ & $\delta$
& $\rho$ & $\delta$ & $\sigma_{v}$\\
Approximate Density &  &  &  &  &  &  &  &  &  & \\\cline{1-1}%
\multicolumn{1}{c}{} &  &  &  &  &  &  &  &  &  & \\
AUKF & \multicolumn{1}{|c}{0.0263} & \multicolumn{1}{c}{0.0535} & 0.0935 &
0.0529 & 0.0370 & 0.0185 & 0.0798 & 0.0201 & 0.0862 & 0.0287\\
Euler & \multicolumn{1}{|l}{0.0096} & 0.0434 & 0.0664 & 0.0072 & 0.0308 &
0.0175 & 0.0818 & 0.0120 & 0.0459 & 0.0242\\
ABC-Joint Score & \multicolumn{1}{|c}{-} & \multicolumn{1}{c}{-} & - &
0.0054 & \textbf{0.0217} & 0.0101 & 0.0441 & \textbf{0.0063} & \textbf{0.0124}
& \textbf{0.0166}\\
ABC-Marginal Score & \multicolumn{1}{|c}{\textbf{0.0028}} &
\multicolumn{1}{c}{0.0301} & \textbf{0.0392} & \textbf{0.0045} & 0.0219 &
\textbf{0.0048} & 0.0480 & 0.0085 & 0.0381 & 0.0167\\
ABC-SS & \multicolumn{1}{|c}{0.0353} & \multicolumn{1}{c}{0.0316} & 0.0427 &
0.0310 & 0.0234 & 0.0119 & \textbf{0.0312} & 0.0109 & 0.0389 & 0.0170\\
ABC-FP & \multicolumn{1}{|c}{0.0178} & \multicolumn{1}{c}{\textbf{0.0215}} &
0.0431 & 0.0145 & 0.0233 & 0.0093 & 0.0358 & 0.0124 & 0.0407 & 0.0168\\
& \multicolumn{1}{|l}{} &  &  &  &  &  &  &  &  & \\\hline\hline
\end{tabular}
%

%TCIMACRO{\TeXButton{E}{\end{table}}}%
%BeginExpansion
\end{table}%
%EndExpansion
%

%TCIMACRO{\TeXButton{end_landscape}{\end{landscape}}}%
%BeginExpansion
\end{landscape}%
%EndExpansion

In Panels B and C respectively of Table 1, we record all RMSE results for the
case when two, then all three parameters are unknown, with a view to gauging
the relative performance of the ABC methods when multiple matches are
required. In these multiple parameter cases, a preliminary run of ABC, based
on uniform priors defined over the domains defined above, has been used to
determine the high mass region of the joint posterior. This information has
been used to further truncate the priors in a subsequent ABC run, and final
marginal posterior estimates then produced. The results recorded in Panels B
and C highlight that when two parameters are unknown (either $\rho$ and
$\sigma_{v}$ or $\rho$ and $\delta$), the score-based ABC method produces the
most accurate density estimates in three of the four cases. Marginalization
produces an improvement in accuracy for $\rho.$ For the other two parameters
marginalization does not yield an increase in accuracy; however, the
differences between the joint and marginal score estimates are minimal. Only
in one case (as pertains to $\delta$) does an ABC method based on summary
statistics outperform the score-based methods. In all four cases, the AUKF
estimate is inferior to all other comparators, and the Euler approximation
also inferior to all ABC-based estimates in three of the four cases. As is
seen in Panel D, when all three parameters are to be estimated, the
score-based ABC estimates remain the most accurate, with the joint score
method superior overall and yielding notably improvements in accuracy over
both the AUKF and Euler approximations.

\section{Conclusions and Discussion\label{end}}

This paper has explored the application of approximate Bayesian computation in
the state space setting. Certain fundamental results have been established,
namely the lack of reduction to finite sample sufficiency and the Bayesian
consistency of the auxiliary model-based method. The (limiting) equivalence of
ABC estimates produced by the use of both the maximum likelihood and
score-based summary statistics has also been demonstrated. The idea of
tackling the dimensionality issue that plagues the application of ABC in high
dimensional problems via an integrated likelihood approach has been proposed.
The approach has been shown to work extremely well in the case in which the
auxiliary model is exact, and to yield some benefits otherwise. However, a
much more comprehensive analysis of different non-linear settings (and
auxiliary models) would be required for a definitive conclusion to be drawn
about the trade-off between the gain to be had from marginalization and the
loss that may stem from integrating over an \textit{inaccurate} auxiliary model.

Indeed, the most important challenge that remains, as is common to the related
frequentist techniques of indirect inference and efficient methods of moments,
is the specification of a computationally efficient and accurate approximating
model. Given the additional need for parsimony, in order to minimize the
number of statistics used in the matching exercise, the principle of aiming
for a large nesting model, with a view to attaining full asymptotic
sufficiency, is not an attractive one. We have illustrated the use of one
simple approximation approach based on the unscented Kalman filter. The
relative success of this approach in the particular example considered,
certainly in comparison with methods based on other more \textit{ad hoc}
choices of summary statistics, augers well for the success of score-based
methods in the non-linear setting. Further exploration of approximation
methods in other non-linear state space models is the subject of on-going
research. (See also Creel and Kristensen, 2014, for some contributions on this front.)

Finally, we note that despite the focus of this paper being on inference about
the static parameters in the state space model, there is nothing to preclude
marginal inference on the states being conducted, at a second stage.
Specifically, conditional on the (accepted) draws used to estimate
$p(\mathbf{\phi|y})$, existing filtering and smoothing methods (including the
recent methods that exploit ABC at the filtering/smoothing level; see, for
example, Jasra \emph{et al., }2010, Calvet and Czellar, 2014, Martin \emph{et
al., }2014) could be used to yield draws of the states, and (marginal)
smoothed posteriors for the states produced via the usual averaging arguments.
With the asymptotic properties of both approaches established (under relevant
conditions), of particular interest would be a comparison of both the finite
sample accuracy and computational burden of the ABC-PMCMC method developed
Martin \emph{et al. }(2014), with that of the method proposed herein, in which
$p(\mathbf{\phi|y})$ is targeted more directly via the score-based approach.

\bigskip

\paragraph{Appendix: Implementation details for the AUKF\newline}

\bigskip

Given the assumed invariance (over time) of both $\nu_{t}$ and $e_{t}$ in
(\ref{discrete_meas}) and (\ref{discrete_state}), the sigma points are
determined as:%
\[
e^{1}=E(e_{t});\text{ }e^{2}=E(e_{t})+a_{e}\sqrt{var(e_{t})};\text{ }%
e^{3}=E(e_{t})-b_{e}\sqrt{var(e_{t})}%
\]
and
\[
v^{1}=E(v_{t});\text{ }v^{2}=E(v_{t})+a_{v}\sqrt{var(v_{t})};\text{ }%
v^{3}=E(v_{t})-b_{v}\sqrt{var(v_{t})}%
\]
respectively, and propagated at each $t$ through the relevant non-linear
transformations, $h_{t}(.)$ and $k_{t}(.).$ The values $a_{e}$, $b_{e}$,
$a_{v}$ and $b_{v}$ are chosen according to the assumed distribution of $e_{t}
$ and $v_{t}$, with a Gaussian assumption for both variables yielding values
of $a_{e}=b_{e}=a_{v}=b_{v}=\sqrt{3}$ as being `optimal'. Different choices of
these values are used to reflect higher-order distributional information and
thereby improve the accuracy with which the mean and variance of the
non-linear transformations are estimated; see Julier \emph{et al. }(2000) and
Ponomareva and Date (2010) for more details. Restricted supports are also
managed via appropriate truncation of the sigma points. The same principles
are applied to produce the mean and variance of the time varying state $x_{t}%
$, except that the sigma points need to be recalculated at each time $t$ to
reflect the up-dated mean and variance of $x_{t}$ as each new value of $y_{t}$
is realized.

In summary, the steps of the AUKF applied to evaluate the likelihood function
of (\ref{discrete_meas}) and (\ref{discrete_state}) are as follows:

\begin{enumerate}
\item Use the (assumed) marginal mean and variance of $x_{t}$, along with the
invariant mean and variance of $v_{t}$ and $e_{t}$ respectively, to create the
$(3\times7)$ matrix of augmented sigma points for $t=0$, $X_{a0}$, as follows.
Define:
\begin{equation}
E(X_{a0})=\left[
\begin{array}
[c]{c}%
E(x_{t})\\
E(v_{t})\\
E(e_{t})
\end{array}
\right]  \text{, }P_{a0}=\left[
\begin{array}
[c]{ccc}%
var(x_{t}) & 0 & 0\\
0 & var(v_{t}) & 0\\
0 & 0 & var(e_{t})
\end{array}
\right]  ,\text{ }\label{e_var}%
\end{equation}
and $\sqrt{P_{a0}}_{j}$ as the $jth$ column of the Cholesky decomposition
(say) of $P_{a0}.$ Given the diagonal form of $P_{a0}$ (in this case), we have%
\[
\sqrt{P_{a0}}_{1}=\left[
\begin{array}
[c]{c}%
\sqrt{var(x_{t})}\\
0\\
0
\end{array}
\right]  ;\text{ }\sqrt{P_{a0}}_{2}=\left[
\begin{array}
[c]{c}%
0\\
\sqrt{var(v_{t})}\\
0
\end{array}
\right]  ;\text{ }\sqrt{P_{a0}}_{1}=\left[
\begin{array}
[c]{c}%
0\\
0\\
\sqrt{var(e_{t})}%
\end{array}
\right]  .
\]
The seven columns of $X_{a0}$ are then generated by%
\[
E(X_{a0});\text{ }E(X_{a0})+a_{j}\sqrt{P_{a0}}_{j}\text{ ; for }j=1,2,3;\text{
}E(X_{a0})-b_{j}\sqrt{P_{a0}}_{j}\text{ ; for }j=1,2,3,
\]
where $a_{1}=a_{x}$, $a_{2}=a_{v}$ and $a_{3}=a_{e}$, and the corresponding
notation is used for $b_{j}$, $j=1,2,3.$

\item Propagate the $t=0$ sigma points through the transition equation as
$X_{x1}=k_{1}\left(  X_{a0},\mathbf{\phi}\right)  $ and estimate the
predictive mean and variance of $x_{1}$ as:%
\begin{align}
E(x_{1}|y_{0})  & =%
%TCIMACRO{\tsum \limits_{i=1}^{7}}%
%BeginExpansion
{\textstyle\sum\limits_{i=1}^{7}}
%EndExpansion
w_{i}X_{x1}^{i}\label{pred_e}\\
var(x_{1}|y_{0})  & =%
%TCIMACRO{\tsum \limits_{i=1}^{7}}%
%BeginExpansion
{\textstyle\sum\limits_{i=1}^{7}}
%EndExpansion
w_{i}(X_{x1}^{i}-E(x_{1}|y_{0}))^{2},\label{pred_var}%
\end{align}
where $X_{x1}^{i}$ denotes the $ith$ element of the $(1\times7)$ vector
$X_{x1}$ and $w_{i}$ the associated weight, determined as an appropriate
function of the $a_{j}$ and $b_{j};$ see Ponomareva and Date (2010).

\item Produce a new matrix of sigma points, $X_{a1},$ for $t=1$ generated by%
\[
E(X_{a1});\text{ }E(X_{a1})+a_{j}\sqrt{P_{a1}}_{j}\text{ ; for }j=1,2,3;\text{
}E(X_{a1})-b_{j}\sqrt{P_{a1}}_{j}\text{ ; for }j=1,2,3,
\]
using the updated formulae for the mean and variance of $x_{t}$ from
(\ref{pred_e}) and (\ref{pred_var}) respectively, in the calculation of
$E(X_{a1}) $ and $P_{a1}$.

\item Propagate the $t=1$ sigma points through the measurement equation as
$X_{y1}=h_{1}\left(  X_{a1},\mathbf{\phi}\right)  $ and estimate the
predictive mean and variance of $y_{1}$ as:%
\begin{align}
E(y_{1}|y_{0})  & =%
%TCIMACRO{\tsum \limits_{i=1}^{7}}%
%BeginExpansion
{\textstyle\sum\limits_{i=1}^{7}}
%EndExpansion
w_{i}X_{y1}^{i}\label{e_y}\\
var(y_{1}|y_{0})  & =%
%TCIMACRO{\tsum \limits_{i=1}^{7}}%
%BeginExpansion
{\textstyle\sum\limits_{i=1}^{7}}
%EndExpansion
w_{i}(X_{y1}^{i}-E(y_{1}|y_{0}))^{2},\label{var_y}%
\end{align}
where $X_{y1}^{i}$ denotes the $ith$ element of the $(1\times7)$ vector
$X_{y1}$ and $w_{i}$ is as defined in Step 3.

\item Estimate the first component of the likelihood function, $p(y_{1}%
|y_{0})$, as a Gaussian distribution with mean and variance as given in
(\ref{e_y}) and (\ref{var_y}) respectively.

\item Given observation $y_{1}$ produce the up-dated filtered mean and
variance of $x_{t}$ via the usual KF up-dating equations:%
\begin{align*}
E(x_{1}|y_{1})  & =E(x_{1}|y_{0})+M_{1}(y_{1}-E(y_{1}|y_{0}))\\
var(x_{1}|y_{1})  & =var(x_{1}|y_{0})-M_{1}^{2}var(y_{1}|y_{0}),
\end{align*}
where:%
\[
M_{1}=\frac{%
%TCIMACRO{\tsum \limits_{i=1}^{7}}%
%BeginExpansion
{\textstyle\sum\limits_{i=1}^{7}}
%EndExpansion
w_{i}(X_{x1}^{i}-E(x_{1}|y_{0}))(X_{y1}^{i}-E(y_{1}|y_{0}))}{var(y_{1}|y_{0})}%
\]
and the $X_{x1}^{i}$, $i=1,2,...,7$ are as computed in Step 3.

\item Continue as for Steps 2 to 6, with the obvious up-dating of the time
periods and the associated indexing of the random variables and sigma points,
and with the likelihood function in evaluated as the product of the components
produced in each implementation of Step 5, and the log-likelihood in
(\ref{approx_like}) produced accordingly.
\end{enumerate}

\end{document}